\documentclass[preprint,12pt]{elsarticle}

\usepackage{amssymb}
\usepackage{amsmath}

\newtheorem{theorem}{Theorem}[section]
\newtheorem{lemma}[theorem]{Lemma}
\newtheorem{corollary}[theorem]{Corollary}
\newtheorem{proposition}[theorem]{Proposition}
\newtheorem{remark}[theorem]{Remark}
\newtheorem{example}[theorem]{Example}
\newenvironment{proof}{\noindent{\em Proof:}}{\quad \hfill$\Box$\vspace{2ex}}
\usepackage{url}
\def \bR {{\mathbb R}}
\def \cP {{\mathcal P}}
\def \be {{\boldsymbol e}}
\def \bb {{\boldsymbol b}}
\def \bx {{\boldsymbol x}}
\def \bz {{\boldsymbol z}}
\def \by {{\boldsymbol y}}
\def \bu {{\boldsymbol u}}

\def \bA {\boldsymbol A}
\def \bI {\boldsymbol I}
\def \bP {\boldsymbol P}

\def \supp {\,{\rm supp}\,}

\def \sign {\,{\rm sign}\,}
\def \prox {\,{\rm prox}\,}
\usepackage{color}
\usepackage{soul}
\usepackage{float}
\usepackage{subcaption}

\begin{document}

\begin{frontmatter}



\title{Computing the proximal operator {\color{black}of the $q$-th power} of the $\ell_{1,q}$-norm for group sparsity}

\author[label1]{Rongrong Lin}
\author[label1]{Shihai Chen}
\author[label2]{Han Feng}
\author[label1]{Yulan Liu\footnote{Corresponding author. \\ E-mail addresses: linrr@gdut.edu.cn (R. Lin), wyyxcsh@163.com (S. Chen), hanfeng@cityu.edu.hk (H. Feng), ylliu@gdut.edu.cn (Y. Liu).}}
\affiliation[label1]{organization={Guangdong University of Technology},
            addressline={Yinglong Road},
            city={Guangzhou},
            postcode={510520},
            country={P.R. China}}

\affiliation[label2]{organization={City University of Hong Kong},
            addressline={Tat Chee Avenue, Kowloon Tong},
            city={Hongkong},
            country={P.R. China}}



\begin{abstract}
In this note, we comprehensively characterize the proximal operator {\color{black}of the $q$-th power} of the $\ell_{1,q}$-norm {\color{black} (denoted by $\ell_{1,q}^{q}$) } with $0\!<\!q\!<\!1$ by exploiting the well-known proximal operator of $|\cdot|^q$ on the real line. In particular, much more explicit characterizations can be obtained whenever $q\!=\!1/2$ and $q\!=\!2/3$ due to the existence of closed-form expressions for the proximal operators of $|\cdot|^{1/2}$ and $|\cdot|^{2/3}$. Numerical experiments demonstrate potential advantages of the $\ell_{1,q}^{ {\color{black} q} }$ regularization in the {\color{black}inter-group and intra-group} sparse vector recovery.
\end{abstract}




\end{frontmatter}


\section{Introduction}
The group-sparsity setting, where the components of the decision vector are grouped together into several distinguishable index sets, has been extensively studied in recent years because of its practical applications in a wide range of fields, including gene selection and analysis \cite{Meier2008},
 cell fate conversion \cite{Hu2024}, signal and image processing \cite{Sun2014,Wang2017}, matrix completion \cite{Liu2023,ZhangLi2024}, multi-task feature learning \cite{LiuJi2009}, multiple kernel learning \cite{Bach2008} and neural networks \cite{Wen2016} in machine learning.
 Many results and notions known for the standard sparse model were generalized to the
block-sparse settings. Among them are the well-known restricted isometry property \cite{Baldassarre2016,Eldar2009}, block convex relaxation techniques \cite{Eldar2009,Eldar2010}, and several $\ell_0$-based
methods (the $\ell_0$-norm counts the number of nonzero blocks) such as the generalized compressive sampling matching pursuit and block iterative hard thresholding  in \cite{Baldassarre2016,Beck2019,Duarte2009}. 
Let $\mathcal{G}:=\{\mathcal{G}_1,\mathcal{G}_2,\ldots,\mathcal{G}_r\}$ be a prescribed non-overlapping partition of the index set $[l]\!:=\!\{1,2,\ldots,l\}$, namely,
$\mathcal{G}_i\cap \mathcal{G}_k\!=\!\emptyset$ for any $i\!\neq\! k\!\in\! [r]$ and $\cup_{i=1}^r \mathcal{G}_i\!=\![l]$.
For any  $\bx\!\in\! \mathbb{R}^l$ and $0\!\leq q\leq 1\leq\! p<+\infty$, let {\color{black}
\(
    \|\bx\|_{p,q}^q\!:=\!\sum_{i=1}^r\|{\bx_{\mathcal{G}_i}}\|_p^q
\)
be the $\ell_{p,q}^q$-norm associated to the partition $\mathcal{G}$, where $\|\bx_{\mathcal{G}_i}\|^q_p\!:=\!(\sum_{j=1}^{|\mathcal{G}_i|} |({\bx}_{\mathcal{G}_i})_j|^p)^{q/p}$ is the $\ell^q_p$-norm of the subvector $\bx_{\mathcal{G}_i}$ with $|\mathcal{G}_i|$ being the cardinality of $\mathcal{G}_i$.
} The group sparse optimization models via the $\ell^q_{p,q}$ regularization in \cite{Liu2023,Hu2017,Pan2021,Yuan2006,Zhang2022} are generally modeled by {\color{black}adding the $\ell_{p,q}^{q}$-norm as the regularized term as follows}:
\begin{align}\label{ProbReg}
    \min\limits_{\bx\in \mathbb{R}^l}\frac{1}{2}\|\bA\bx-\bb\|^2_2+\lambda \|\bx\|_{p,q}^q,
\end{align}
where the measurement matrix $\bA\!\in\!\bR^{m\times l}$ and the observation vector $\bb\in\bR^m$ are given and $\lambda\!>\!0$ is the regularization parameter.
Particularly, the model \eqref{ProbReg} with  $p\!=\!2$ and $q\!=\!1$  is known as  group Lasso problem,
firstly introduced for grouped variable selection in statistics \cite{Yuan2006}.
Under some group restricted eigenvalue conditions, 
the global and local recovery bounds for \eqref{ProbReg}  were developed in \cite{Hu2017}. Later, a general lower recovery bound for \eqref{ProbReg} with $p\!=\!2$ was established in \cite{Feng2020}. 
Recently, $\ell_{1,0}$-norm regularized optimization was built  for multivariate regression \cite{Chen2023}. They provided a thorough analysis of the connection between the minimized  $\ell_{1,0}$-norm problem, the $\ell_{1,0}$ regularized optimization and the corresponding  problem  with $\ell_{1,0}$ constraints  regarding stationary points, local minimizers and global minimizers.

In order to solve the problem \eqref{ProbReg}, numerous effective methods have been presented, and the proximal gradient approach is one of the most prevalent optimization algorithms.
Proximal algorithms are fundamental tools to handle optimization problems involving non-smooth functions or large-scale datasets \cite{Beck2017,Parikh2014}. The core of proximal algorithms lies in the evaluation of the proximal operator of a function. Given a function $f:\bR^n\to\bR$ and a parameter $\nu>0$, the proximal operator $\prox_{\nu f}(\by)$ at a point $\by\in\bR^n$ is a set-valued mapping defined by
\begin{align}\label{ProxDef}
\prox_{\nu f}(\by):=\arg\min_{\bu\in\bR^n} f(\bu)+\frac{1}{2\nu}\|\bu-\by\|_2^2.
\end{align}
An increasing number of researchers are concentrating on characterizing the proximal operator of sparsity-inducing functions such as  
the $\ell_0$-norm \cite{Blumensath2008}, the $\ell_1$-norm \cite{Beck2009}, the $\ell^q_q$-norm with $0<q<1$ \cite{Cao2013,Chen2016}, Capped $\ell_1$ \cite{ZhangTong2010}, Transformed $\ell_1$-norm \cite{Zhang2017}, SCAD \cite{Fan2001}, MCP \cite{Zhang2010}, PiE \cite{Liu2024}, log-sum \cite{Prater2022}, the difference of the $\ell_1$-norm and $\ell_2$-norm \cite{Lou2018}, the ratio $\ell_1/\ell_2$ of the $\ell_1$-norm and $\ell_2$-norm \cite{Tao2022}, $(\ell_1/\ell_2)^2$ \cite{Jia2024}, the $p$-power $\|\cdot\|_1^p$ of $\ell_1$-norm with $p>1$ \cite{Prater2023}, and so on. Specially, for the optimization problem \eqref{ProbReg}, 
the proximal operators of the $\ell_{2,q}^{q }$-norm with some particular values $q\!\in\!\{\frac12,\frac23\}$ have been characterized in \cite[Proposition 18 (iii) and (v)]{Hu2017}.
For a general value $0\!<\!q\!<\!1$, the proximal operators of the $\ell_{2,q}^{{\color{black} q }}$-norm were derived in {\textcolor{black}{\cite[Proposition 4.2]{Hu2024}}}  and \cite[Remark 3.7]{Lin2025} in two different ways. The latter is a more compact characterization by directly revealing the close relationship between the proximal operators of the $\ell_{2,q}^{{\color{black} q }}$-norm and the $\ell_q^q$-norm.
The proximal operators of the $\ell_{1,q}^{{\color{black} q }}$-norm  with some particular values $q\!\in\!\{\frac12,\frac23\}$ have been listed without proof in \cite[Proposition 18 (iv) and (vi)]{Hu2017}. Unfortunately, these results are incorrect. Two counterexamples are constructed, see Example \ref{counterexample12} and Example \ref{Exmp2} in Section \ref{Section3}. 

In this paper, we 
 focus on the computation and understanding of the {\textcolor{black}{proximal}} operator of {\color{black} $\ell_{1,q}^{q}$} with $0\!<\!q\!<\!1$. In Section \ref{Section2}, some preliminaries including the existing characterizations of the proximal operator of the $\ell_q$-norm and some necessary properties of the proximal operator of $\|\cdot\|_1^{q}$ are given.
 In Section \ref{Section3},  the expression of $\prox_{\nu\|\cdot\|_1^{q}}$
 is derived by exploiting the well-known proximal operator of the $\ell^q_q$-norm on $\bR$.
Some further properties of $\prox_{\nu\|\cdot\|_1^{q}}$ are developed as products.
Thanks to all cubic polynomials and quartic polynomials are solvable by radicals \cite[p. 73-75]{Zwillinger2018},  more explicit expressions for $\prox_{\nu\|\cdot\|_1^{q}}$ can be obtained whenever $q\!=\!1/2$ and $q\!=\!2/3$, respectively. In Section \ref{Section4}, we conduct group sparse vector recovery based on the $\ell_{1,q}^{{\color{black} q}}$ regularization.
Conclusions and remarks are given in the final section.

 \section{Preliminaries}\label{Section2}
 
Throughout the paper, write $\bR^n_{\downarrow}:=\{\bz\!\in\!\bR^n:z_1\ge z_2\ge \cdots\ge z_n\ge0\}$. For any $\bz\in\bR^n_{\downarrow}$ and $\kappa\in[n]$, denote by $\|\bz\|_{(\kappa)}:=\sum_{i=1}^{\kappa}|z_i|$ the Ky-Fan $\kappa$-norm of $\bz$ \cite[Example 7.24]{Beck2017}.
For any $\bz\in \mathbb{R}^n$ and  positive constants $v>0$,  denote by $\|\bz\|_v:=(\sum\limits_{i=1}^{n} |z_i|^v)^{1/v}$ the $\ell_v$-norm and denote by $\supp(\bz):=\{i\in[n]: z_i\ne0\}$ the support of $\bz$.
Given {\textcolor{black}{an}} index subset $T$ of $[n]$, denote by $\bz_T\in\bR^{|T|}$ the subvector of $\bz$ by deleting the components that are not in $T$ from $\bz$, where $|T|$ stands for the cardinality of the set $T$.
Let $\cP_n$ denote the set of all $n\times n$ signed permutation matrices, that is, those matrices that have only one nonzero entry in every row or column, which is $\pm 1$.  For any  $\by,\bz\in\bR^n$, the pareto inequaltiy $\by\succ \bz$ (respectively, $\by\!\succeq\! \bz$) means that $y_i\!>\! z_i$  (respectively, $y_i\!\ge\! z_i$) for all $i\!\in\![n]$.  Let $\by\!\nsucc\! \bz$ denote $\by\!\succ\! \bz$ is not satisfied.
Define the signum function $\sign(\bz)$ on $\bR^n$ as  
\begin{equation*}
[\sign(\bz)]_i:=\left\{\begin{array}{ll}
1,& \text{ if } z_i>0, \\
0,  &\text{ if }  z_i=0,\\
-1,& \text{ otherwise },\\
\end{array}\right.\ i\in[n],
\end{equation*} 
and ${\rm Sign}(\bz)$ on $\bR^n_{}$ as  \begin{equation*}
{\color{black}{[{\rm Sign}(\bz)]_i}}:=\left\{\begin{array}{ll}
\{1\},& \text{ if } z_i>0, \\
{\rm [-1,1]},  & \text{ if } z_i=0,\\
\{-1\},& \text{ otherwise},\\
\end{array}\right.\ i\in[n].
\end{equation*} 
Denote by ${\bf 0}_n$ and ${\bf 1}_n$ the  vectors in $\bR^n$ with all components being $0$ and $1$, respectively.
For any $i\in[n]$, denote by $\be_i\in\bR^n$  an unit vector with the $i$-th coordinate being $1$.

  \subsection{The proximal operator of the $\ell_q^q$-norm}
  
 {\color{black}Recall that the regularization term $\|\bx\|^q_{p,q}\!=\!\sum_{i=1}^r \|\bx_{\mathcal{G}_i}\|_1^q$ for any $\bx\in\bR^l$ in \eqref{ProbReg} is separable. In particular, when $r=l$, the $\ell_{1,q}^q$-norm 
  simplifies to the $\ell_{q}^q$-norm.
 By the definition of the proximal operator and Theorem 6.6 in \cite{Beck2017}, it suffices to characterize the proximal operator $\prox_{\nu \|\cdot\|_1^q}$ of $\|\cdot\|_1^q$ on $\bR^{n}$ with a parameter $\nu>0$ and $n\ge2$. Specifically, we have
 $$
 \prox_{\nu\|\cdot\|^q_{1,q}}(\bx)= \prox_{\nu\|\cdot\|_1^q}(\bx_{\mathcal{G}_1})\times \prox_{\nu\|\cdot\|_1^q}(\bx_{\mathcal{G}_2})\times\cdots\times \prox_{\nu\|\cdot\|_1^q}(\bx_{\mathcal{G}_r}).
 $$
 It will be observed in Proposition \ref{Case2Prop} that $\prox_{\nu \|\cdot\|_1^q}$ has a close connection with the proximal operator of the $\ell_{q}^{q}$-norm on $\bR$, namely, $|\cdot|^q$. For this reason, we need to recall the characterization of the proximal operator of $|\cdot|^q$.}
The proximal operator of $|\cdot|^q$ with a general value $0<q<1$ has been fully derived in \cite{Chen2016}. 
For completeness, we rewrite it as follows.
\begin{lemma}\label{LpProx} Let $0\!<\!q\!<\!1$ and $\nu\!>\!0$. It holds that for any  $\tau\!\in\!\mathbb{R}$,
\begin{align*}
    \prox_{\nu |\cdot|^{q}}(\tau)=
    \begin{cases}
    \{0\}, &\text{ if }|\tau|<c_{\nu,q},\\
    \{0,\sign(\tau)\rho_{\nu,q}\}, &\text{ if }|\tau|=c_{\nu,q},\\
    \{\sign{(\tau)}r(\tau)\}, & \text{ otherwise},
    \end{cases}
\end{align*}
where $c_{\nu,q}\!=\!\frac{2-q}{2(1-q)}(2\nu(1-q))^{1/(2-q)}$, $\rho_{\nu,q}\!=\!(2\nu(1-q))^{1/(2-q)}$, $r(\tau)\!\in\! ((\nu q(1-q))^{1/(2-q)},|\tau|)$ is the unique root of the equation
\begin{equation}\label{eqlpprox}
r(\tau)+\nu q r(\tau)^{q-1}-\tau=0.
\end{equation}
In particular, $\prox_{\nu |\cdot|^{q}}$ has a closed-form expression whenever $q\!=\!1/2$ and $q\!=\!2/3$ as follows:
$$
\prox_{\nu|\cdot|^{1/2}}(\tau)=\left\{\begin{array}{ll}
    \{0\}, & \mbox{ if }  |\tau|<\frac{3}{2}\nu^{\frac{2}{3}},\\
    \{0,\sign(\tau) \nu^{\frac{2}{3}}\}, & \mbox{ if }  |\tau|=\frac{3}{2}\nu^{\frac{2}{3}}, \\
   \Big\{\frac{4\tau}{3}\cos^2\Big(\frac{1}{3}\arccos\big(-\frac{3\sqrt{3}}{4}\nu|\tau|^{-\frac{3}{2}})\Big)\Big\},  &  \text{ otherwise},\\
\end{array}
\right.
$$
and 
$$
\prox_{\nu|\cdot|^{2/3}}(\tau)=\left\{\begin{array}{ll}
    \{0\}, &  \mbox{ if } |\tau|<2(\frac{2}{3}\nu)^{\frac{3}{4}},\\
    \{0,\sign(\tau)(\frac{2}{3}\nu)^{\frac{3}{4}}\}, &   \mbox{ if }|\tau|=2(\frac{2}{3}\nu)^{\frac{3}{4}}, \\
   \Big\{\frac{1}{8}\sign(\tau)\big(\sqrt{2t}+\sqrt{\frac{2|\tau|}{\sqrt{2t}}-2t}\big)^3\Big\},  &  \text{ otherwise},\\
\end{array}
\right.
$$
with $t=\Big(\frac{\tau^2}{16}+\sqrt{\frac{\tau^4}{256}-\frac{8\nu^3}{729}}\Big)^{1/3}+\Big(\frac{\tau^2}{16}
-\sqrt{\frac{\tau^4}{256}-\frac{8\nu^3}{729}}\Big)^{1/3}$.
\end{lemma}

We remark that the proximal operator $\prox_{\nu |\cdot|^{q}}$ in Lemma \ref{LpProx} for an arbitrary value $p\notin\{\frac12,\frac23\}$ can be numerically evaluated by the Newton method \cite[Algorithm 1]{Chen2016} or the bisection method \cite[Algorithm 1]{Liu2024b}. In addition, from the definition of the  proximal operator  \eqref{ProxDef} and Lemma \ref{LpProx}, the following result is easily checked.
\begin{corollary} \label{lqProp}
  Let $0\!<\!q\!<\!1$ and $\nu\!>\!0$. it holds that
  $\max \prox_{\nu|\cdot|^{q}}(\tau^1)<\prox_{\nu|\cdot|^{q}}(\tau^2)$ for any  $\tau^2>\tau^1\geq c_{\nu,q}$ and 
   \(
  \prox_{\nu\alpha|\cdot|^{q}}(\alpha\tau)=\alpha\prox_{\nu\alpha^{q-1}|\cdot|^{q}}(\tau)
  \)
   for any $\alpha>0$ and $\tau\in \mathbb{R}$.
\end{corollary}

\subsection{$\prox_{\nu\|\cdot\|_1^q}$'s general properties }

In the subsection, we assume that $0\!<\!q\!<\!1$ and $\nu\!>\!0$. 
 By the definition of proximal operator \eqref{ProxDef}, the proximal operator $\prox_{\nu \|\cdot\|_1^q}$ of the $\ell_{1}^{q}$-norm  with a parameter $\nu$  follows:
\begin{align}\label{Proxl1qDef}
\prox_{\nu \|\cdot\|_1^q}(\by)=\arg\min_{\bu\in\bR^n} J_{\by}(\bu):=\nu\|\bu\|_1^q+\frac12\|\bu-\by\|_2^2, \,\mbox{ for any } \by\in \mathbb{R}^n.
\end{align}
By abuse of notation, if the set $\prox_{\nu\|\cdot\|_1^{q}}(\by)$ has a single element $\tilde{\by}$, we denote $\prox_{\nu\|\cdot\|_1^{q}}(\by)=\{\tilde{\by}\}$ or simply  $\prox_{\nu\|\cdot\|_1^{q}}(\by)=\tilde{\by}$. Applying the first optimal condition \cite[Theorem 10.1]{RW09} to the problem \eqref{Proxl1qDef}, if ${\bf 0}_n\ne\tilde{\by}\in \prox_{\nu\|\cdot\|_1^q}(\by)$, then $\tilde{\by}\in S_q(\by)$, where
\begin{align}\label{SsetDef}
S_q(\by):=\{\bu\in \mathbb{R}^n\,:\, \by\in \bu+\nu q\|\bu\|_1^{q-1} {\rm Sign}(\bu)\}.
\end{align}
 Obviously, ${\bf 0}_n\!\notin S_q(\by)\!$ for any $\by\ne {\bf 0}_n$ since $0\!<\!q\!<\!1$.  
Hence, 
\begin{align}\label{ProxEq}
\prox_{\nu\|\cdot\|_1^q}(\by)=\arg\min_{\bu\in\{{\bf 0}_n\}\cup S_q(\by)}J_{\by}(\bu).
\end{align}
The following assertions directly follow from the definition of $S_q(\by)$ in \eqref{SsetDef}.

 \begin{proposition}\label{GeneralProp}
 Given ${\bf 0}_n\!\neq\! \by\!\in\! \mathbb{R}^n$. Then, for  any $\bu\!\in\! S_q(\by)$, defined as in \eqref{SsetDef}, the following assertions hold.
 \begin{itemize}
     \item [(i)] 
     One has $u_i=y_i-\nu q \|\bu\|_1^{q-1} \sign(u_i)$ for any  $i\in \supp(\bu)$. Moreover,  $\sign(u_i)=\sign(y_i)$ and $|u_i|<|y_i|$ for any $i\in {\rm supp} (\bu)$.
     
     \item [(ii)] For any $i\in [n]$,  the following implication holds:
     \begin{align*}
         |y_i|\leq \nu q\|\bu\|_1^{q-1}\Longrightarrow u_i=0.
     \end{align*}
  \end{itemize}
  \end{proposition}
 
 \begin{proposition} \label{Case3Prop}
 Given $\by\in\bR^n$. If $\|\by\|_{\infty}:=\max\limits_{i\in [n]}|y_i|>c_{\nu,q}$, then ${\bf 0}_n\notin\prox_{\nu\|\cdot\|_1^q}(\by)$,  where $c_{\nu,q}$ is defined as in Lemma \ref{LpProx}. 
\end{proposition}
\begin{proof} Assume that $\|\by\|_{\infty}>c_{\nu,q}$. Then, there exists an index $i_0\in[n]$ such that $|y_{i_0}|>c_{\nu,q}$.  By  Lemma \ref{LpProx}, the function $\nu|t|^q+\frac12(t-y_{i_0})^2$ has the unique  and nonzero minimizer, namely,
\begin{equation}\label{ui0}
\min_{t\in\bR}\nu |t|^q+\frac12(t-y_{i_0})^2< \frac12y_{i_0}^2.
\end{equation}
Then, by the definition of $J_{\by}(\bu)$ in \eqref{Proxl1qDef}, we compute
$$
\begin{array}{ll}
\displaystyle{\min_{\bu\in\bR^n}J_{\by}(\bu) }
& \displaystyle{ =\min_{\bu\in\bR^n}\nu\|\bu\|_1^q+\frac12\|\bu-\by\|_2^2} \displaystyle{\le \min_{\bu\in\bR^n}\nu\|\bu\|_q^q+\frac12\|\bu-\by\|_2^2 }\\
&\displaystyle{=\min_{u_{i_0}\in\bR} \nu|u_{i_0}|^q+(u_{i_0}-y_{i_0})^2+\sum_{i_0\ne j\in[n]}\min_{u_j\in\bR}\nu|u_j|^q+(u_j-y_j)^2 }\\
&\displaystyle{< \sum_{i=1}^n\frac12y_i^2 =\frac12\|\by\|_2^2=J_{\by}({\bf 0}_n),}
\end{array}
$$
where the first inequality is by $\|\bu\|_1\le\|\bu\|_q$ for any $\bu\in\bR^n$, and the second inequality is obtained by \eqref{ui0} and the fact that $\min\limits_{t\in\bR}\nu |t|^q+\frac12(t-y_{j})^2\le \frac12y_j^2$ for any $i_0\ne j\in[n]$. Thus, the minimizer of $J_{\by}{\color{black} (\bu) }$ must be nonzero, which completes the proof.
\end{proof}

\begin{proposition}\label{lemmapermutation} Let $\bP\in\cP_n$. Then $\prox_{\nu\|\cdot\|_1^q}(\by)\!=\! \bP^{-1}\prox_{\nu\|\cdot\|_1^q}(\bP\by)$ for any $\by\in\bR^n$, where $\bP^{-1}$ denotes the inverse matrix of $\bP$. 
\end{proposition} 
\begin{proof} Notice that $\|\cdot\|_1^q$ is a signed permutation invariant function on $\bR^n$. The desired result  can be derived directly from  the definition of $\prox_{\nu\|\cdot\|_1^q}$ in \eqref{Proxl1qDef}.
\end{proof}

By Proposition \ref{lemmapermutation}, it suffices to consider $\by\in\bR^n_{\downarrow}$. From now on, we assume  $\by\in \bR^n_{\downarrow}$.

\begin{proposition}\label{partialorder} 
Given $\by\in \bR^n_{\downarrow}$. If
${\bf 0}_n\notin \prox_{\nu\|\cdot\|_1^q}(\by)$, then ${\bf 0}_n\notin \prox_{\nu\|\cdot\|_1^q}(\bar{\by})$ for any $\bar{\by}\succ \by$. 
\end{proposition}
\begin{proof} Take any $\bar{\by}\succ \by$. Then there exist $\bz\succ {\bf 0}_n$ such that $\bar{\by}=\by+\bz$. 
By Equation \eqref{Proxl1qDef}, we have
$$
\begin{array}{ll}
\displaystyle{ \min_{\bu\in\bR^n}J_{\bar{\by}}(\bu)  }
&\displaystyle{ =\min_{\bu\in\bR^n}\nu\|\bu\|_1^q+\frac12\|\bu-(\by+\bz)\|_2^2 } \le J_{\bar{\by}}(\tilde{\by}) \\
&=\nu\|\tilde{\by}\|_1^q+\frac12\|(\by-\tilde{\by})+\bz\|_2^2 \\
&=\nu\|\tilde{\by}\|_1^q+\frac12\|\by-\tilde{\by}\|_2^2  +\frac12\|\bz\|_2^2+\bz^{\top}(\by-\tilde{\by})\\
&=J_{\by}(\tilde{\by})+\frac12\|\bz\|_2^2+\bz^{\top}(\by-\tilde{\by})\\
&\leq J_{\by}({\bf 0}_n)+\frac12\|\bz\|_2^2+\bz^{\top}(\by-\tilde{\by})\\
&= \frac12\|\by\|_2^2 + \frac12\|\bz\|_2^2+\bz^{\top}\by-\bz^{\top}\tilde{\by}\\ &=\frac12\|\bar{\by}\|_2^2-\bz^{\top}\tilde{\by} \\
&<\frac12\|\bar{\by}\|_2^2=J_{\bar{\by}}({\bf 0}_n),
\end{array}
$$
where the first inequality is obtained by taking $\tilde{\by}\!\in\! \prox_{\nu\|\cdot\|_1^q}(\by)$ with $\tilde{\by}\neq {\bf 0}_n$, and the last inequality is by $\bz^{\top}\tilde{\by}>0$ as $\bz\succ {\bf 0}_n$ and ${\bf 0}_n\ne\tilde{\by}\succeq {\bf 0}_n$.
\end{proof}
\begin{remark} Fix $n\!=\!2$, $\nu\!=\!1$, and $\by\!\in\![0,+\infty)^2:=[0,+\infty)\times[0,+\infty)$. For $q=1/2$ and $q=2/3$,  Figures \ref{ProxL11over2} (c) and  \ref{ProxL12over3} (c) coincide with the results of Propositions  \ref{Case3Prop}  and \ref{partialorder}.  In Figure \ref{ProxL11over2} (c), for instance, $\prox_{\|\cdot\|_1^{1/2}}(\by)\!=\!\{{\bf 0}_2\}$ whenever $\by$ is in the blue region. For any vector $\by$ outside the closure of the blue region,  which is a subset of the square $[0,c_{\nu,q}]^2=[0,\frac{3}{2}]^2$ from Proposition \ref{Case3Prop}, it holds that ${\bf 0}_2\notin\prox_{\|\cdot\|_1^{1/2}}(\by)$ by Proposition \ref{partialorder}.
\end{remark}

\subsection{ $\prox_{\nu\|\cdot\|_1^q}$'s properties when sparsity is given}

In this subsection, we focus on 
 analyzing the property of $\tilde{\by}\in \prox_{\nu\|\cdot\|_1^q}(\by)$ for any $\by\in \bR^n_{\downarrow}$  when the sparsity is given and provide an explicit formula to characterize the values of $a:=\|\tilde{\by}\|_1$ when $\|\tilde{\by}\|_0$ is given. 
 First, we need some important lemmas.
 Given $0<s\leq n$, write
\begin{align}\label{t0Def}
   \widehat{t}_{\nu,q}(s):=(\nu s q (1-q))^{1/(2-q)},\quad
    \widetilde{t}_{\nu,q}(s):=(2-q)\Big(\frac{\nu s q}{(1-q)^{1-q}}\Big)^{1/(2-q)}.
 \end{align}
 For simplicity, we adopt $\widehat{t}(s)$ and $\widetilde{t}(s)$ instead of $\widehat{t}_{\nu,q}(s)$ and $\widetilde{t}_{\nu,q}(s)$ unless otherwise stated.
 
  \begin{lemma}\label{aExsitLemma} Fix $\nu\!>\!0$ and $0\!<\!q\!<\!1$.
  Given $0\!<\!s\!\leq\!n$ and $b\!\geq\! 0$. Suppose that $ g(t)\!:=\!t+\nu sq t^{q-1}-b$ 
   and  $h(t)\!:=\!\nu t^{q}+\frac{t^2}{2s}-\frac{b}{s}t$ 
   for any $t\in \mathbb{R}$. 
   \begin{itemize}
       \item [(i)]  The function $g$ is strictly  convex on $(0,\infty)$ and achieves {\color{black} its} unique  global minimum value at $t\!=\!\widehat{t}(s)$ with minimal value  $g(\widehat{t}(s))\!=\!\widetilde{t}(s)-b$ where 
    $\widehat{t}(s)$ and $\widetilde{t}(s)$ are defined as in \eqref{t0Def}. In addition, if $b\!\geq\! \widetilde{t}(s)$, then $g(t)\!=\!0$ has at most two roots, assumed to be $t_1$ and $t_2$ respectively,  satisfying $0\!<\!t_1\!\leq\! \widehat{t}(s)\!\leq\! t_2<b$.
    
    \item [(ii)] If  $b\geq \widetilde{t}(s)$, then $h$ is strictly increasing on the intervals $[0,t_1]$ and $[t_2,+\infty)$ and strictly decreasing on $[t_1,t_2]$. Therefore, $h(t_2)\leq h(t_1)$.
   \end{itemize}
  \end{lemma}
 \begin{proof}
 Notice that $g'(t)=1+\nu sq(q-1)t^{q-2}$
  and $g''(t)=\nu sq(q-1)(q-2)t^{q-3}$.
  Obviously, the conclusion (i) holds.
  After calculation, we know that
  $$
  h'(t)=\frac{\nu s q t ^{q-1}+t-{\color{black} b}}{s}=\frac{g(t)}{s}\quad {\rm and }\quad  h''(t)=\frac{g'(t)}{s}.
  $$
  Hence, (ii) follows from (i).
 \end{proof}
 
 By Proposition \ref{GeneralProp} and Lemma \ref{aExsitLemma}, we can derive the following conclusions.
 
 \begin{lemma}\label{SpLemma}
 Given ${\color{black} {\bf 0}_n}\!\neq\! \by\!\in\! \mathbb{R}^n$. Then, for  any $\bu\in S_q(\by)$ with $s:=\|\bu\|_0$ and $a:=\|\bu\|_1$, the following assertions hold.
 \begin{itemize}
     \item [(i)]  For any $i,j\in [n]$ and  $u_j>0$, if $y_i>y_j$, then 
      $u_i>u_j$ and if $y_i=y_j$, then 
      $u_i=u_j$. In a word, if $\by\in \mathbb{R}^n_{\downarrow}$, then
      $\bu_{[s]}$ has the same descending order with $\by_{[s]}$.
    
     \item [(ii)] Suppose that $\by\in \mathbb{R}^n_{\downarrow}$. If $\|\by\|_{(s)}\geq \widetilde{t}(s)$, then $a=\|\bu\|_1$ is a solution on $(0, \infty)$ of the equation $\nu q s t^{q-1}+t=\|\by\|_{(s)}$ with
     $ y_{s}>\nu q a^{q-1}$ and $y_i\leq \nu q a^{q-1}$ for any $i=s+1,\ldots,n$.
     Moreover,
   $u$ is characterized by
    \begin{align*}
    u_{i}=\left\{\begin{array}{cl}
     y_{i}-\nu q a^{q-1} , & {\rm  if\;  } i\in [s],\\
       0, &  {\rm otherwise}.
     \end{array}
    \right.
    \end{align*}
    Otherwise, there is no vector $\bu$ with $\|\bu\|_0=s$ in the set $S_q(\by)$ if $\|\by\|_{(s)}< \widetilde{t}(s)$.
  \end{itemize}
  \end{lemma}
 
If $\by$ is a multiple of $\be_i$ or ${\bf 1}_n$, we can explicitly reveal the connection between the proximal operators of the $\ell_{1,q}^{q }$-norm and $|\cdot|^q$ on $\bR$. 

 \begin{corollary}
\label{Lemmabasis}  For any $\tau\!>\!0$, the following statements hold.
\begin{itemize}
 \item[(i)]
 $\prox_{\nu \|\cdot\|_1^{q}}(\tau {\bf 1}_n)\!=\!\{\tau^*{\bf 1}_n\,|\, \tau^*\!\in \!\prox_{\nu n^{q-1} |\cdot|^{q}}(\tau)\}$. 
 \item[(ii)] $\prox_{\nu \|\cdot\|_1^{q}}(\tau \be_i)\!=\!\{\tau^*\be_i\,|\, \tau^*\!\in \!\prox_{\nu |\cdot|^{q}}(\tau)\}$ for any $i\in[n]$.
 \end{itemize}
 \end{corollary}
 \begin{proof} We first prove Statement (i).
 Take any $\tilde{\by}\in \prox_{\nu \|\cdot\|_1^{q}}(\tau {\bf 1}_n)$.
 Activating  Lemma \ref{SpLemma} with $\by=\tau {\bf 1}_n$ yields that $\tilde{\by}_i=\tilde{\by}_j:=\tau_0$ for any $i,j\in [n]$. Hence, $\tilde{\by}=\tau_0 {\bf 1}_n$.
 Notice that
\begin{align*}
\min_{\bu\in\bR^n}\nu\|\bu\|_1^q+\frac12\|\bu-\tau{\bf 1}_n\|_2^2&=\nu n^q |\tau_0|^q+ \frac{1}{2}n(\tau_0-\tau)^2\geq \min_{t\in\bR} \nu n^q |t|^q+ \frac{1}{2}n(t-\tau)^2.
\end{align*}
On the other hand, take any $\tau^*\!\in \!\prox_{\nu n^{q-1} |\cdot|^{q}}(\tau)$ and write $\widehat{\by}:=\tau^* {\bf 1}_n$. Then
\begin{align*}
  \min_{t\in\bR} \nu n^q |t|^q+ \frac{1}{2}n(t-\tau)^2&=\nu n^q |\tau^*|^q+ \frac{1}{2}n(\tau^*-\tau)^2=\nu\|\widehat{\by}\|_1^q+\frac12\|\widehat{\by}-\tau{\bf 1}_n\|_2^2\\
  &\geq \min_{\bu\in\bR^n}\nu\|\bu\|_1^q+\frac12\|\bu-\tau{\bf 1}_n\|_2^2.
\end{align*}
The above two inequalities yield that
$\min\limits_{\bu\in\bR^n}\nu\|\bu\|_1^q+\frac12\|\bu-\tau{\bf 1}_n\|_2^2=\min\limits_{t\in\bR} \nu n^q |t|^q+ \frac{1}{2}n(t-\tau)^2$. 

The proof of Statement (ii) can be done in a similar way. 
\end{proof}

\begin{proposition}\label{aCharacProp}
Let $\by\!\in\!\bR^n_{\downarrow}$. If $\tilde{\by}\!\in\!\prox_{\nu\|\cdot\|_1^q}(\by)$ with $s\!:=\!\|\tilde{\by}\|_0\!\neq\! 0$, then 
$\|\tilde{\by}\|_1\!=\!\prox_{\nu s|\cdot|^q}(\|\by\|_{(s)})$.
\end{proposition}
\begin{proof}
 Take any ${\bf 0}_n\ne\tilde{\by}\in\prox_{\nu\|\cdot\|_1^q}(\by)$. Then $\tilde{\by}\in S_q(\by)$.
  Let $a:=\|\tilde{\by}\|_1$. From Lemma \ref{SpLemma}, it follows that $\|\by\|_{(s)}\geq \widetilde{t}(s)$, $\tilde{\by}_{[s]}=\by_{[s]}-\nu qa^{q-1}{\bf 1}_s$ and  
 \begin{equation}\label{Equ1}
 \nu q s a^{q-1}+a-\|\by\|_{(s)}=0.
 \end{equation}
 Write $g(t)\!:=\! \nu q s t^{q-1}+t-\|\by\|_{(s)}$ and 
 $h(t)\!:=\!\nu t^{q}+\frac{t^2}{2s}-\frac{\|\by\|_{(s)}}{s}t$ for any $t\!>\!0$.
Clearly, $0\!<\!a<\!\|\by\|_{(s)}$  and $a$ is a root of $g(t)\!=\!0$. Let $\overline{[s]}:=\{s+1,s+2,\dots,n\}$ be the complement of $[s]$ in $[n]$.
On the other hand, by $\tilde{\by}_{[s]}=\by_{[s]}-\nu qa^{q-1}{\bf 1}_s$ and \eqref{Equ1}, we have 
\begin{equation*}
\begin{array}{ll}
\min\limits_{\bu\in\bR^n} J_{\by}(\bu)=J_{\by}(\tilde{\by})
&=\nu\|\tilde{\by}\|_1^q+\frac{1}{2}\|\tilde{\by}-\by\|_2^2 \\
&=\nu\|\tilde{\by}\|_1^q+\frac12\|\nu q\|\tilde{\by}\|_1^{q-1}\sign(\by_{[s]})\|_2^2+\frac12\|\by_{\overline{[s]}}\|_2^2\\
&=\nu a^q+\frac12 (\nu q)^2 sa^{2q-2}+\frac12\|\by_{\overline{[s]}}\|_2^2\\
&=\nu a^q+\frac{(\|\by\|_{(s)}-a)^2}{2s}+\frac12\|\by_{\overline{[s]}}\|_2^2=h(a)+\frac{\|\by\|_{(s)}^2}{2s}+\frac12\|\by_{\overline{[s]}}\|_2^2,
\end{array}
\end{equation*}
which implies that $\widehat{t}(s)\leq a\leq \|\by\|_{(s)}$ by activating Lemma \ref{aExsitLemma} with $b=\|\by\|_{(s)}$.
Hence, $a$ is the unique root of $g(t)=0$.  Replacing $r$, $\nu$, and $\tau$  in \eqref{eqlpprox} of Lemma \ref{LpProx} by $a$, $\nu s$, and $\|\by\|_{(s)}$, respectively, we obtain
$a=\prox_{\nu s|\cdot|^q}(\|\by\|_{(s)})$.
The desired result is obtained.
\end{proof}

 For any $\by\in \mathbb{R}^n_{\downarrow}$ and $s\in [n]${\color{black},} define
 \begin{equation}\label{Spys}
S_q(\by,s)\!:=\!\! \left\{\begin{array}{ll}
\{\max\{\by\!-\!c_s{\bf 1}_n,{\bf 0}_n\}\}, &\mbox{ if }  \|\by\|_{(s)}\!\ge\! \widetilde{t}(s), y_s\!>\!c_s, y_{i}\!\le\! c_s, \forall i\in \{s+1,\ldots,n\},\\
\emptyset & {\rm otherwise},
\end{array}
  \right.\end{equation}
 where $\widetilde{t}(s)$ is defined as in \eqref{t0Def} and  if $\|\by\|_{(s)}\!\ge\! \widetilde{t}(s)$, $c_s$ depending on $s$ is given by
\begin{equation}\label{Spycs}
 c_s:=\nu q  (\prox_{\nu s|\cdot|^q}(\|\by\|_{(s)}))^{q-1}.
 \end{equation}

\begin{proposition}\label{Case2Prop}  For any $\by\!\succ\! c_{\nu,q}{\bf 1}_n$, we have $\prox_{\nu\|\cdot\|_1^q}(\by)\!=\!\by-\nu q  (\prox_{\nu n|\cdot|^q}(\|\by\|_1))^{q-1}{\bf 1}_n$, where $c_{\nu,q}$ is defined as in Lemma \ref{LpProx}.
\end{proposition}
\begin{proof} 
Fix any $\by\succ c_{\nu,q}{\bf 1}_n$. By Proposition \ref{Case3Prop}, ${\bf 0}_n\notin\prox_{\nu\|\cdot\|_1^q}(\by)$. Take any $\tilde{\by}\in\prox_{\nu\|\cdot\|_1^q}(\by)$ with
$\widetilde{s}:=\|\tilde{\by}\|_0$. The desired result can be deduced if $\widetilde{s}=n$ by Proposition \ref{aCharacProp} and Lemma \ref{SpLemma}(ii).
In fact, conversely, suppose that there exist $\overline{\bu}\in \prox_{\nu\|\cdot\|_1^q}(\by)$ and $\overline{s}=\|\overline{\bu}\|_0<n$.
Then $\overline{\bu}\in S_q(\by,\overline{s})$ by Proposition \ref{aCharacProp} and Lemma \ref{SpLemma}(ii). Hence, $\|\by\|_{(\overline{s})}\!\ge\! \widetilde{t}(\overline{s})$,  $y_{\overline{s}}>c_{\overline{s}}$
and $y_i\leq c_{\overline{s}}$ for any $i=\overline{s}+1,\ldots,n$
by \eqref{Spys}. On the other hand,  notice that
$\|\by\|_{(\overline{s})}>\overline{s}c_{\nu,q}$ since $ \by\succ c_{\nu,q}{\bf 1}_n$. 
So,  by  Corollary \ref{lqProp} and Lemma \ref{LpProx}, we have that
\[\prox_{\nu {\overline{s}}|\cdot|^q}(\|\by\|_{(\overline{s})})> \max\prox_{\nu \overline{s}|\cdot|^q}(\overline{s}c_{\nu,q})= \overline{s}\max\prox_{\nu \overline{s}^{q-1}|\cdot|^q}(c_{\nu,q})=\overline{s} \rho_{\nu\overline{s}^{q-1},q}=\rho_{\nu\overline{s},q}
\]
where $\rho_{\kappa,q}\!:=\!(2\kappa(1-q))^{1/(2-q)}$.
Associating this with
\eqref{Spycs} yields that
\(
c_{\overline{s}}<\nu q(\rho_{\overline{s}\nu,q})^{q-1}.
\)
Hence, summarily, for any $i=\overline{s}+1,\ldots,n$, it holds that
$c_{\nu,q}<y_i\leq c_{\overline{s}}<\nu q(\rho_{\overline{s}\nu,q})^{q-1}$,
 that is,
$$
\nu q (2\nu \overline{s} (1-q))^{\frac{q-1}{2-q}}   >\frac{2-q}{2(1-q)}(2\nu(1-q))^{\frac{1}{2-q}}.
$$
After simple calculation, the above inequality is 
 reduced to $q \overline{s}^{-\frac{1-q}{2-q}}>2-q$,
 which contradicts the inequality $q {\color{black}\overline{s}}^{-\frac{1-q}{2-q}}\le q<2-q$ for any $0<q<1$. The proof is complete.
\end{proof}

\section{The proximal operator of the $\ell_{1,q}^{{\color{black} q }}$-norm}\label{Section3}

Many  algorithms have been proposed to solve the problem \eqref{ProbReg} {\color{black}with $p=1$}, and one of the most popular optimization algorithms is the proximal gradient method, whose $k+1$ iteration step is
\begin{align*}
    \bx^{k+1}:={\rm arg}\min\limits_{\bx\in \mathbb{R}^n} \big\{\lambda\|\bx\|_{1,q}^q+\frac{1}{2\nu}\|\bx-\by^k\|_2^2\big\}
    {\rm\ \  with\ \ }
    \by^k=\bx^k-\nu \bA^{\top} (\bA\bx^k-\bb),
\end{align*}
 where $\nu>0$ is a parameter.
Associating this with \eqref{ProxDef} yields that
\[
\bx^{k+1}=\prox_{\lambda\nu \|\cdot\|_1^q}(\by^{k}_{\mathcal{G}_1})\times \prox_{\lambda\nu \|\cdot\|_1^q}(\by^{k}_{\mathcal{G}_2})\times\cdots\times \prox_{\lambda\nu \|\cdot\|_1^q}(\by^{k}_{\mathcal{G}_r}).
\]
For any $k$ and any $i\in [r]$, there exists $\bP^k_{i}\in \mathcal{P}_{|\mathcal{G}_i|}$ such
that $\by^k_i:=\bP^k_{i}\by^{k}_{\mathcal{G}_i}\in \mathbb{R}^{|\mathcal{G}_i|}_{\downarrow}$. Associating the above equality with {\color{black} Proposition} \ref{lemmapermutation} yields that
\begin{align}\label{MainProx}
\bx^{k+1}\!=\!(\bP^k_{1})^{-1}\prox_{\lambda\nu \|\cdot\|_1^q}(\by^{k}_1)\times (\bP^k_{2})^{-1}\prox_{\lambda\nu \|\cdot\|_1^q}(\by^{k}_2)\times\cdots\times (\bP^k_{r})^{-1}\prox_{\lambda\nu \|\cdot\|_1^q}(\by^{k}_r).
\end{align}
Hence, to compute $\bx^{k+1}$, the characterization of $\prox_{\lambda\nu \|\cdot\|_1^q}(\by^k_i)$ is crucial for any $\by^k_i\in \mathbb{R}^{|\mathcal{G}_i|}_{\downarrow}$.

\subsection{Main results}

 Given a point $\by\in\bR^n_{\downarrow}$,
 from the definition of $S_q (\by)$ in \eqref{SsetDef}, $S_q (\by,s)$ in \eqref{Spys} and properties discussed in Lemma \ref{SpLemma}, it follows that
 \begin{equation}\label{Spy}
 S_{\color{black} q }(\by)=\bigcup_{s\in [n]} S_{\color{black} q }(\by,s). 
 \end{equation}
Recall 
$$
\prox_{\nu\|\cdot\|_1^q}(\by)=\arg\min_{\bu\in\{{\bf 0}_n\}\cup S_q (\by)}J_{\by}(\bu)
$$
 in \eqref{ProxEq}. 
The following characterization follows.

\begin{theorem}\label{MainTh} 
Given a point ${\bf 0}_n \neq \by\in\bR^n_{\downarrow}$. It holds that
$$
\prox_{\nu\|\cdot\|_1^q}(\by)=\left\{\begin{array}{ll}
\{{\bf 0}_n\}, & \mbox{ if } \|\by\|_1\le  \frac{2-q}{1-q}(\nu q(1-q))^{\frac{1}{2-q}},\\
{\color{black}\{} \by-c_n {\bf 1}_n {\color{black}\}}, & \mbox{ if } \by\succ c_{\nu,q}{\bf 1}_n,\\
\arg\min\limits_{\bu\in S_q(\by)}J_{\by}(\bu), & \mbox{ if } \|\by\|_{\infty}> c_{\nu,q}\mbox{ and } \by\nsucc c_{\nu,q}{\bf 1}_n,\\
\arg\min\limits_{\bu\in\{{\bf 0}_n\}\cup S_q(\by)}J_{\by}(\bu),&  \mbox{ otherwise},
\end{array}
\right.
$$
where $S_q(\by)$ is given by \eqref{Spy}, $c_{\nu,q}$ is defined as in Lemma \ref{LpProx} and $c_s$ given as in \eqref{Spycs}.
\end{theorem}
\begin{proof} 
If $\|\by\|_1< \frac{2-q}{1-q}(\nu q(1-q))^{\frac{1}{2-q}}$, then  $\widetilde{t}(s)-\|\by\|_{(s)}\geq \widetilde{t}(1)-\|\by\|_1>0$ for any $s\geq 1$, where $\widetilde{t}(s)$ is defined as in \eqref{t0Def}. Hence, $S_q(\by)=\emptyset$ from Lemma \eqref{SpLemma} (ii) and then $\prox_{\nu\|\cdot\|_1^q}(\by)=\{{\bf 0}_n\}$.
 The second and third cases can be deduced by Proposition \ref{Case2Prop} and \ref{Case3Prop}, respectively. 
\end{proof}

 \begin{remark}
  Let $\by\!\in\!\bR^n_{\downarrow}$ with $\|\by\|_{\infty}>c_{\nu,q}$ and write  
  \[
  I_{\color{black} q }(\by):=\{s\in[n] |\,  \|\by\|_{(s)}\!\ge\! \widetilde{t}(s), y_s\!>\!c_s, y_{i}\!\le\! c_s, \forall i\in \{s+1,\ldots,n\}\},
  \]
  and $\overline{s}\!:=\!\max\{s\!\in\!  I_q(\by)\}$. Then $\tilde{\by}:=\max\{\by-c_{\overline{s}}{\bf 1}_n,{\bf 0}_n\}\!\in\!\prox_{\nu\|\cdot\|_1^q}(\by)$.  In fact, 
  by Theorem \ref{MainTh}, we know that ${\bf 0}_n\!\notin\! \prox_{\nu\|\cdot\|_1^q}(\by)$
  and $S_q(\by)\neq \emptyset$. So, 
  $I_q(\by)\neq \emptyset$. If $I_q(\by)$ has a unique element, then the conclusion directly holds.
If $I_q(\by)$ has multiple elements, we take any $s^1< s^2=\overline{s}\in I_q(\by)$ and write  $\bu^i\!:=\!\max\{\by-c_{s^i}{\bf 1}_n,{\bf 0}_n\}$ for $i=1,2$.
Then $\bu^i\in  S_q(\by,s^i)$ for $i=1,2$.
Clearly, $\|\bu^{i}\|_0=s^i$ for $i=1,2$.
Define $\tilde{\bu}^{\epsilon}\in\bR^n$ by
$$
    \bu^{\epsilon}_i:=\left\{\begin{array}{ll}
    \bu^{1}_i,   & 1\le i\le s^1,   \\
    \epsilon,  & s^1<  i\le s^2,\\ 
    0, & \mbox{ otherwise}.
\end{array}
\right. \text{ with any}\; \epsilon>0.
$$
Notice that $\|\bu^{\epsilon}\|_0\!=\!\|\bu^{s^2}\|_0\!=\!s^2$.
From Lemma  \ref{SpLemma}, $\bu^2$ is the unique minimizer of $J_{\by}$ over the set $\{\bz\in\bR^n:\|\bz\|_0\!=\!s^2\}$. Thus,  $J_{\by}(\bu^2)\!<\!J_{\by}(\bu^{\epsilon})$.
By the continuity of $J_{\by}$ and $\lim\limits_{\epsilon\to0}\bu^{\epsilon}\!=\!\bu^1$, we have
\(
J_{\by}(\bu^2)\!\le\! \lim_{\epsilon\to0}J_{\by}(\bu^{\epsilon})\!=\!J_{\by}(\bu^1)\!=\!J_{\by}(\tilde{\by}).
\)
So,  $\min\limits_{\bu\in S_q(\by)}J_{\by}(\bu)\!=\!J_{\by}(\tilde{\by})$.
\end{remark}

\begin{figure}[!htpb]
  \centering
  \begin{subfigure}[b]{2.3in}
    \centering
    \includegraphics[width=2in]{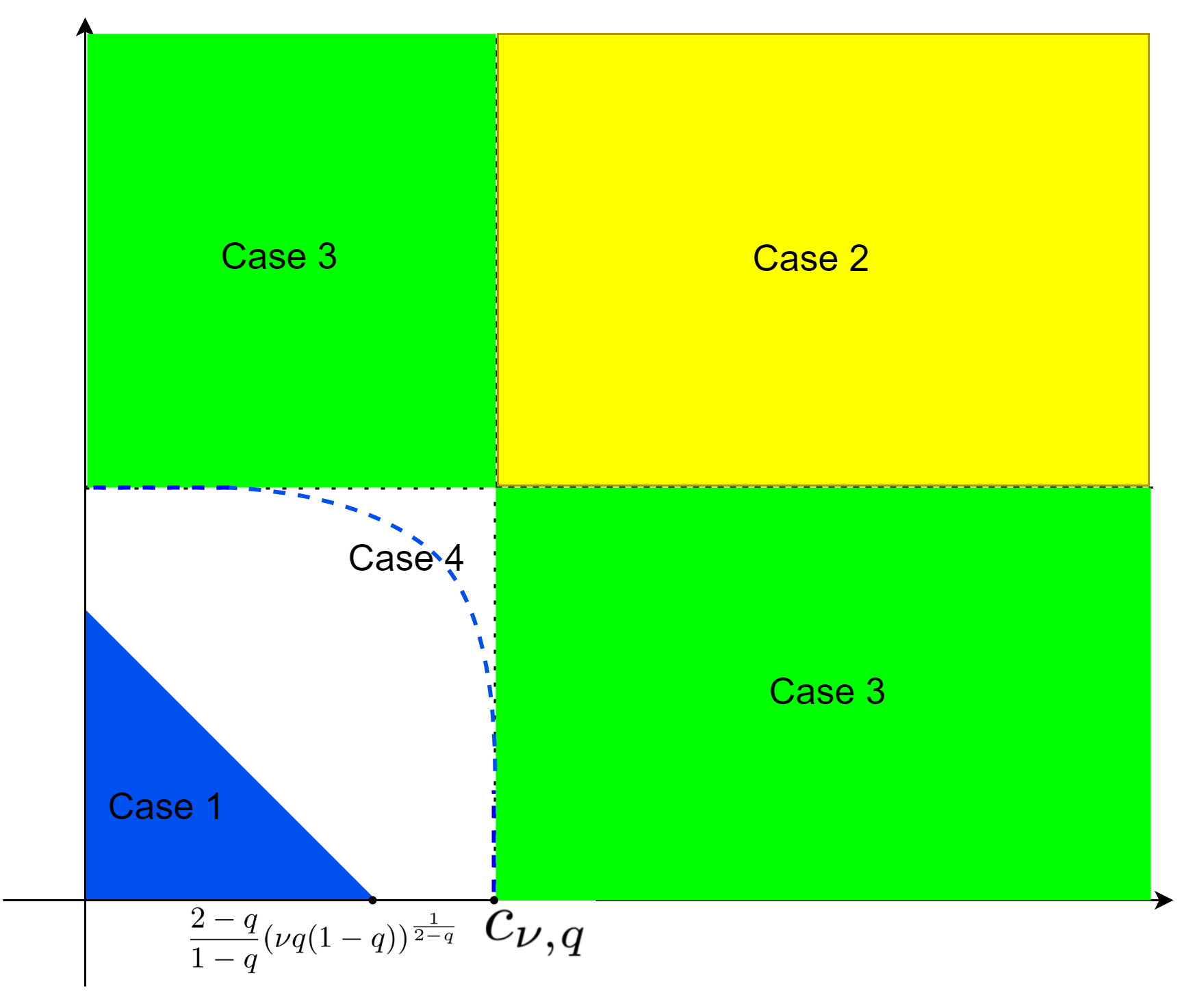}
    \subcaption{Four cases in Theorem \ref{MainTh}}
    \label{fig:sub_a}
  \end{subfigure}%
  \hfill
  \begin{subfigure}[b]{1.6in}
    \centering
    \includegraphics[width=1.7in]{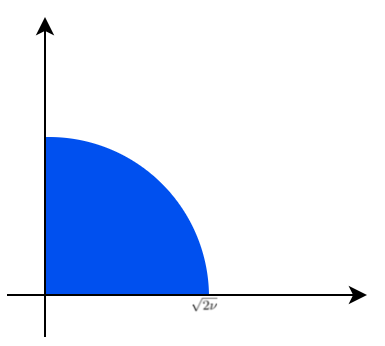}
    \subcaption{$\prox_{\nu\|\cdot\|_1^0}(\by)\!=\!\{{\bf 0}_2\}$}
    \label{fig:sub_b}
  \end{subfigure}
  \hfill
  \begin{subfigure}[b]{1.6in}
    \centering
    \includegraphics[width=1.7in]{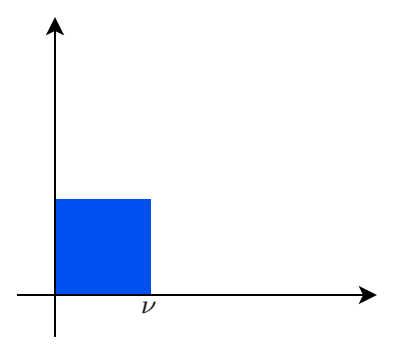}
    \subcaption{$\prox_{\nu\|\cdot\|_1^1}(\by)\!=\!\{{\bf 0}_2\}$}
    \label{fig:sub_c}
  \end{subfigure}
  
  \caption{Illustration of the proximal operator $\prox_{\nu\|\cdot\|_1^q} (\by)$ for any $\by\in[0,+\infty)^2$.}
  \label{Illustration}
\end{figure}

\begin{remark} 
By the definition of  the proximal operator \eqref{ProxDef}, $\prox_{\nu\|\cdot\|_1^0}(\by)$ has the following form:
$$
  \prox_{\nu\|\cdot\|_1^0}(\by)=\left\{\begin{array}{ll}
    \{{\bf 0}_n\}, & \mbox{ if } \|\by\|_2<\sqrt{2\nu},  \\
    \{{\bf 0}_n,\by\}, & \mbox{ if }\|\by\|_2=\sqrt{2\nu},  \\ 
    \{\by\}, & \mbox{ if }\|\by\|_2>\sqrt{2\nu}, \\ 
 \end{array}
 \right.
$$
and 
$
\prox_{\nu\|\cdot\|^1_1}(\by)\!=\!(\sign(y_i)\max\{|y_i|-\nu,0\}: i\in[n])
$
in \cite{Beck2009}. Observe that $\prox_{\nu\|\cdot\|_1^0}(\by)\!=\!\{{\bf 0}_n\}$ whenever $\by\!\in\! \{\bx\in\bR^n : \|\bx\|_2\!<\!\sqrt{2\nu}\}$ and $\prox_{\nu\|\cdot\|^1_1}(\by)\!=\!\{{\bf 0}_n\}$ whenever $\by\!\in\! \{\bx\!\in\!\bR^n \!:\! \|\bx\|_{\infty}\!\le\!\nu\}$, see Figure \ref{Illustration} (b) and (c), respectively. 
 \end{remark}
 
Fix $n\!=\!2$ and $\by\in[0,+\infty)^2$.
The four cases in Theorem \ref{MainTh} are marked in blue, yellow, green and white, and the blue curve is the boundary of the set of points $\by$ such that $\prox_{\nu\|\cdot\|_1^q}(\by)\!=\!\{{\bf 0}_2\}$, as depicted in Figure \ref{Illustration} (a).

{\color{black}
\begin{remark} It is worth mentioning that the proximal operator of 
the $\ell^q_1$-norm with $q>1$ was recently investigated in \cite{Prater2023}. In particular, its proximal operator has a closed-form expression for the power $q\in\{2,3,4\}$. As has been commented in \cite{Prater2023}, $\|\cdot\|_1^q$ with $q>1$ is not a sparsity promoting function. Nevertheless, Theorem \ref{MainTh} delves into the case where $0<q<1$, which is different from the one examined in \cite{Prater2023}. Moreover, our analyses for characterizing the proximal operator of $\|\cdot\|_1^q$ with $0<p<1$ are also significantly different from those in \cite{Prater2023}.
\end{remark}
}

  \subsection{ Characterization of $\prox_{\nu\|\cdot\|_1^{1/2}}$}

  By \eqref{Spy} with $q=1/2$, after simple calculation, we have 
  \begin{equation}\label{S12y}
  S_{1/2}(\by)\!=\!\bigcup_{s\in [n]}\!\big\{\max\{\by-c_s{\bf 1}_n ,{\bf 0}_n \}:  \|\by\|_{(s)}\!\ge\! \big(\frac{3}{2^{\frac43}}\nu^{\frac23}\big)s^{\frac23}, y_s\!>\!c_s,  y_{i}\!\le\! c_s, \forall i\!\in\! \{s+1,\ldots,n\}\big\},
  \end{equation}
  where $c_s$ if $\|\by\|_{(s)}\!\ge\! \Big(\frac{3}{2^{\frac43}}\nu^{\frac23}\Big)s^{\frac23}$ is given by
  \begin{equation}\label{S12ycs}
  c_s \!=\! \frac{\nu}{2} \Big(\prox_{\nu s|\cdot|^{1/2}}(\|\by\|_{(s)})\Big)^{-1/2}\!=\!\frac{\sqrt{3}\nu}{4\sqrt{\|\by\|_{(s)}} \cos\big(\frac{1}{3}\arccos(-\frac{3\sqrt{3}\nu s}{4}\|\by\|_{(s)}^{-\frac32}) \big)}.
  \end{equation}
  Activating  Theorem \ref{MainTh} with $q\!=\!1/2$, we  obtain the following characterization for $\prox_{\nu\|\cdot\|_1^{1/2}}$.
  \begin{corollary}\label{Corollary1over2}
  Given $\by\in \mathbb{R}^n_{\downarrow}$. Then
  \begin{align}\label{q=1/2Resut}
  \prox_{\nu\|\cdot\|_1^{1/2}}(\by)=
  \left\{\begin{array}{ll}
  \{{\bf 0}_n\}, & \mbox{ if } \|\by\|_1\le \frac{3}{2^{\frac43}}\nu^{\frac{2}{3}}, \\
  \{\by-c_n{\bf 1}_n\}, & \mbox{ if }\by\succ \frac{3}{2}\nu^{\frac23}{\bf 1}_n,\\
  \arg\min\limits_{\bu\in S_{1/2}(\by)}J_{\by}(\bu), &\mbox{ if }\|\by\|_{\infty}>\frac{3}{2}\nu^{\frac23}\mbox{ and } \by\nsucc \frac{3}{2}\nu^{\frac23}{\bf 1}_n,\\  
  \arg\min\limits_{\bu\in\{{\bf 0}_n\}\cup S_{1/2}(\by)} J_{\by}(\bu),&\mbox{ otherwise},\\ 
  \end{array}
  \right.
  \end{align}
  where $S_{1/2}(\by)$ is given by \eqref{S12y} and $c_n$ is defined as in \eqref{S12ycs}.
  \end{corollary}
  
\begin{remark}
Hu el. at \cite[Proposition 18 (iii)]{Hu2017} or \cite[Proposition 4.1 (iv)]{Hu2024} proposed the following expression of $\prox_{\nu\|\cdot\|_1^{1/2}}(\by)$ for any $\by\in \mathbb{R}^n$ without proof: 
\begin{equation}\label{wrongL112}
    \prox_{\nu\|\cdot\|_1^{1/2}}(\by)=
    \left\{\begin{array}{ll}
    \{{\bf 0}_n\}, & \mbox{ if } J_{\by}(\tilde{\by}) >J_{\by}({\bf 0}_n),\\
    \{{\bf 0}_n,\tilde{\by}\}, & \mbox{ if }J_{\by}(\tilde{\by})=J_{\by}({\bf 0}_n),\\
    \{\tilde{\by}\}, & \mbox{ if }J_{\by}(\tilde{\by}) <J_{\by}({\bf 0}_n),
    \end{array}
\right.
\end{equation}
where  $\tilde{\by}:=\by-c_n {\rm sign} (\by)$ with $c_n$ is by \eqref{S12ycs}. Notice that 
$\{\tilde{\by}\}=S_{1/2}(\by,n)$, where $S_{1/2}(\by,n)$ is by \eqref{Spys}. 
Hence, 
when $S_{1/2}(\by,n)\!=\! \arg\min\limits_{\bu\in S_{1/2}(\by)} J_{\by}(\bu)$, 
\eqref{q=1/2Resut} is equivalent to \eqref{wrongL112}, however the former is more specific; otherwise, the expression \eqref{wrongL112} is an error, see the following example.
 \end{remark}

\begin{example}\label{counterexample12} (Counterexample) Fix $\nu=1$ and $n=2$. Choose $\by=(\frac{5}{3},\frac{1}{3})^{\top}$.
Take $s=1\in[2]$, then $\|\by\|_{(s)}=|y_1|=\frac{5}{3}\ge \frac{3}{2^{\frac43}}=\widetilde{t}_{1,1/2}(1)\approx 1.19$. By \eqref{S12ycs}, we has
$$
c_1=\frac{\sqrt{3}}{4\sqrt{\frac53} \cos\big(\frac{1}{3}\arccos(-\frac{3\sqrt{3}}{4}(\frac53)^{-\frac32}) \big)}\approx 0.454.
$$
Clearly, $y_1>c_1$ and $y_2\le c_1$. Hence,
$S_{1/2}(\by,1)=\{\bu^1\}$ with
$\bu^1:=(\frac53-c_1,0)^{\top}\approx (1.2126,0)^{\top}$ and  $J_{\by}(\bu^{1})=(1.2126+0)^{\frac12}+\frac12((5/3-1.2126)^2+(1/3)^2)\approx 1.2598$.
Take $s=2\in [2]$, then $\|\by\|_{(s)}=\|\by\|_1=\frac53+\frac13=2>\frac{3}{2^{\frac23}}=\widetilde{t}_{1,1/2}(2)=1.89$. By \eqref{S12ycs}, we arrive at
$$
c_2=\frac{\sqrt{3}}{4\sqrt{2}}\cos\Big(\frac13 \arccos\Big(-\frac{3\sqrt{6}}{8}\Big)\Big)=\frac12.
$$
Since $y_2=\frac{1}{3}< c_2=\frac12$,  $S_{1/2}(y,2)=\emptyset$. Hence, $S_{1/2}(\by)=\{\bu^1\}$ by \eqref{Spy}.
Notice that
$J_{\by}({\bf 0}_2)=\frac12((5/3)^2+(1/3)^2)=\frac{13}{9}\approx 1.444$. Then
 $\prox_{\nu\|\cdot\|_1^{1/2}}(\by)=\{\bu^1\}$ as $J_{\by}(\bu^1)<J_{\by}({\bf 0}_2)$.

However, let $\bu^2:=\by-c_2(1,1)^{\top}=(\frac76,-\frac16)^{\top}$ and then $J_{\by}(\bu^2)=(7/6+1/6)^{\frac12}+\frac12((5/3-7/6)^2+(1/3+1/6)^2)\approx 1.405<J_{\by}({\bf 0}_2)$.
Hence from  \eqref{wrongL112}, it follows that $\prox_{\nu\|\cdot\|_1^{1/2}}(\by)=\{\bu^2\}$, which is obviously incorrect as $J_{\by}(\bu^1)<J_{\by}(\bu^2)$.
\end{example}

\begin{figure}[!htpb]
  \centering
  \begin{subfigure}[b]{0.31\textwidth}
    \centering
    \includegraphics[width=\linewidth]{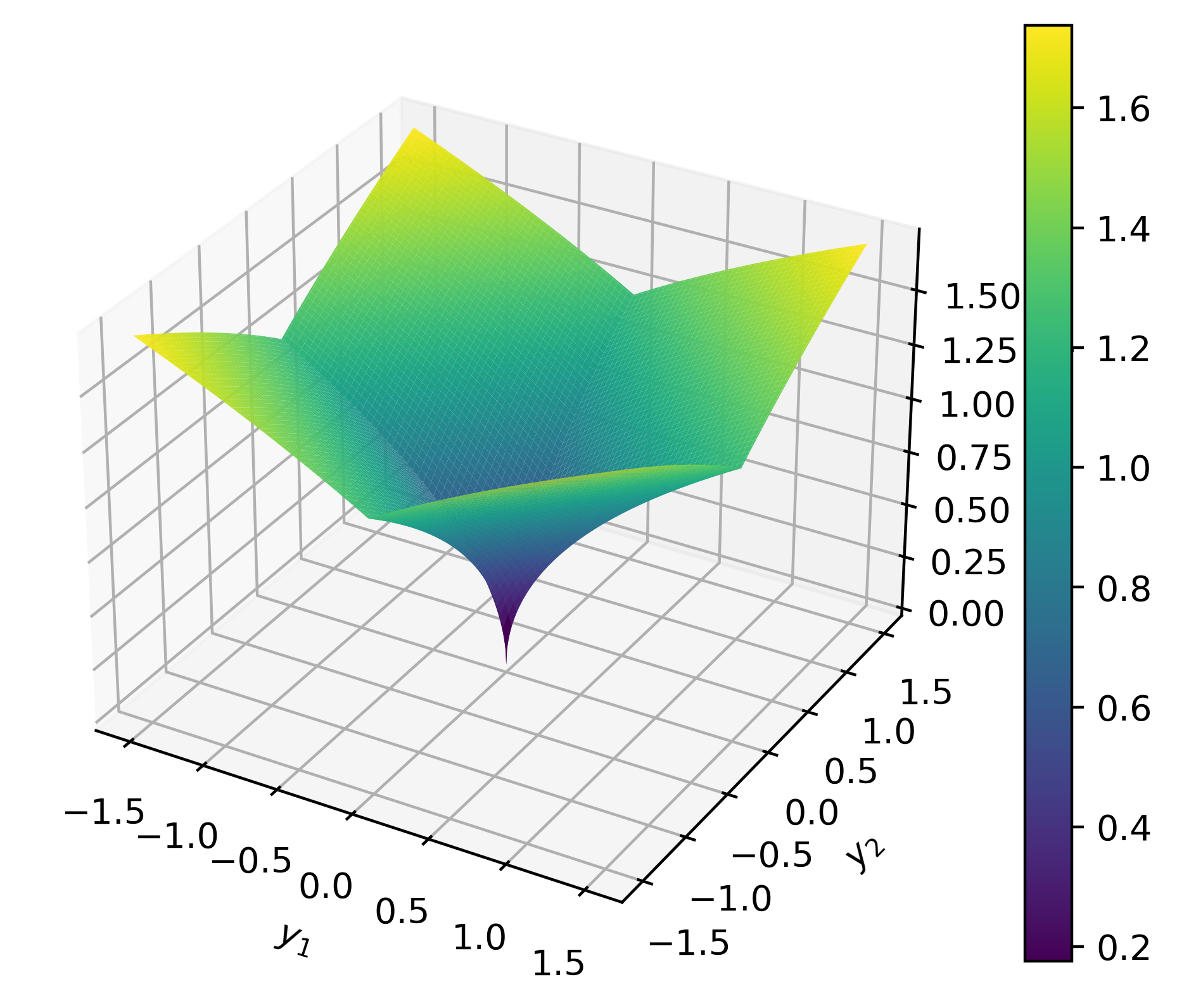}
    \subcaption{Image of $\|\cdot\|_1^{1/2}$ on $\by\in[-2,2]^2$}
    \label{ProxL11over2:sub1}
  \end{subfigure}%
  \hfill
  \begin{subfigure}[b]{0.3\textwidth}
    \centering
    \includegraphics[width=\linewidth]{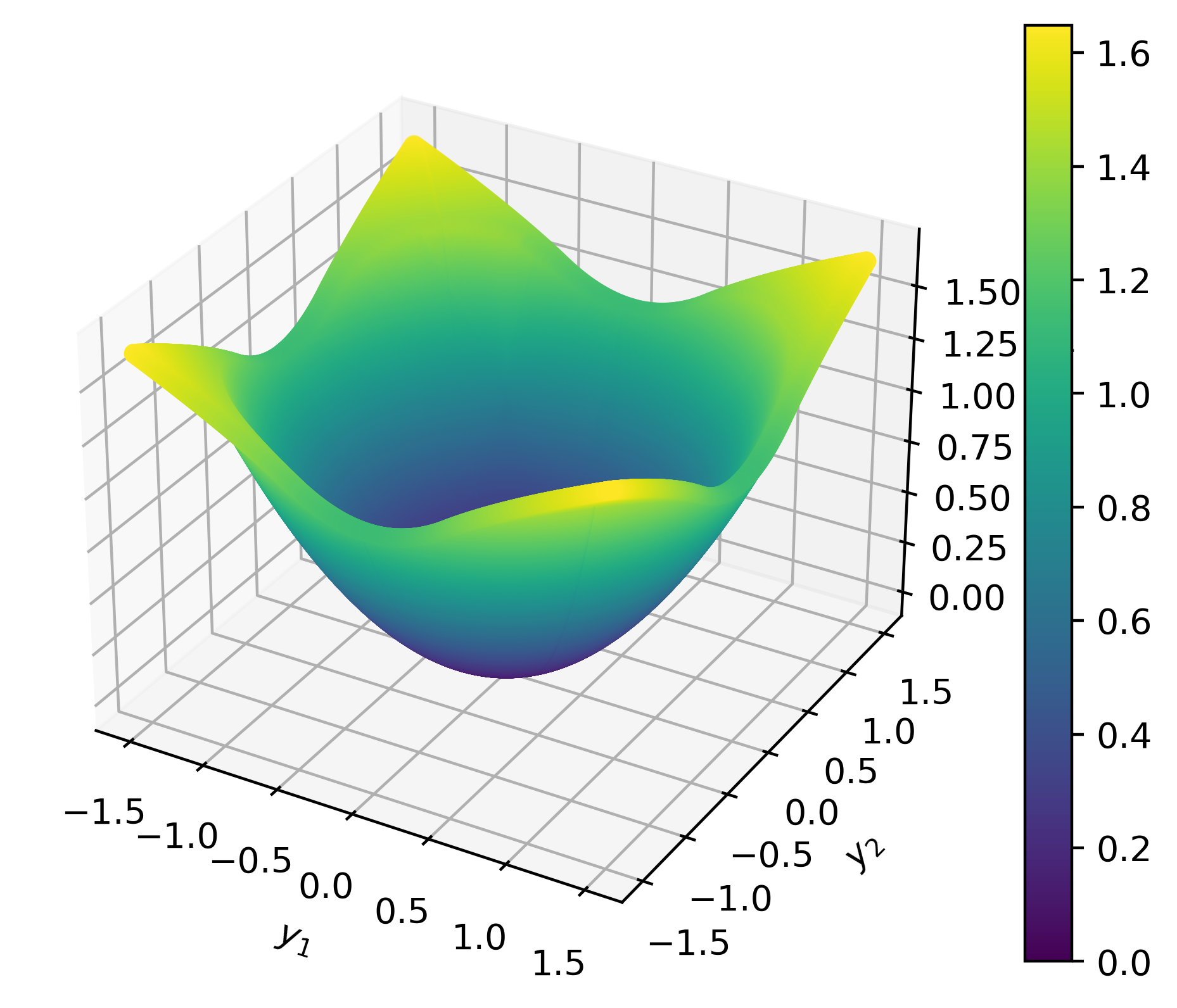}
    \subcaption{Its Moreau envelope}
    \label{ProxL11over2:sub2}
  \end{subfigure}%
  \hfill
  \begin{subfigure}[b]{0.25\textwidth}
    \centering
    \includegraphics[width=\linewidth]{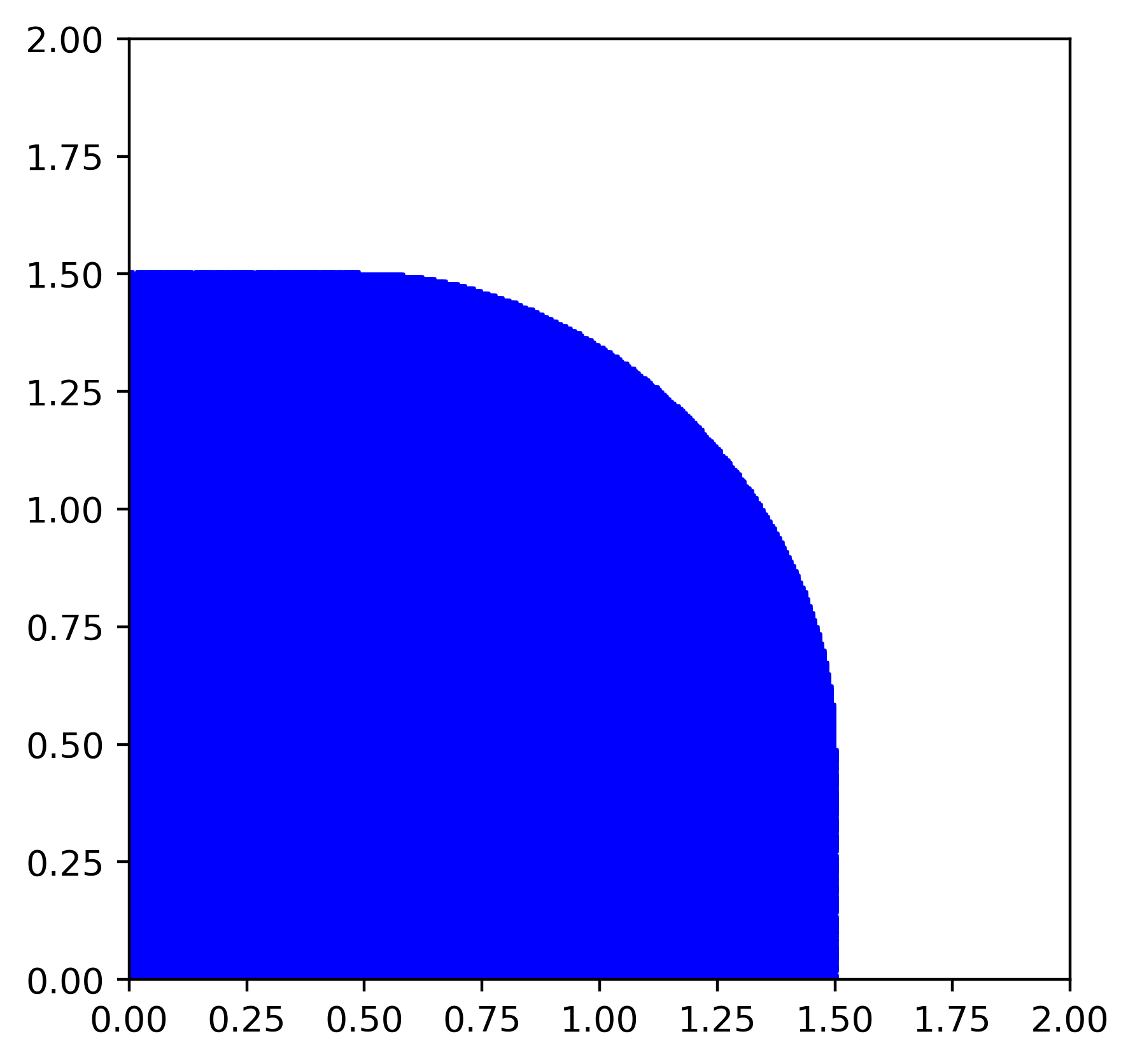}
    \subcaption{$\prox_{\|\cdot\|_1^{1/2}}(\by)\!=\!\{{\bf 0}_2\}$}
    \label{ProxL11over2:sub3}
  \end{subfigure}
  
  \caption{Illustration of $\|\cdot\|_1^{1/2}$, its Moreau envelope and proximal operator.}
  \label{ProxL11over2}
\end{figure}

To further understand the behavior of $\prox_{\nu\|\cdot\|_1^{1/2}}$, we should make some comments. Fix $\nu=1$. From Lemmas \ref{LpProx} and Corollary  \ref{Lemmabasis},
 $\prox_{\|\cdot\|_1^{1/2}}(\tau (1,0)^{\top})=\{{\bf 0}_2\}$ when $\tau<3/2$; and $\prox_{\|\cdot\|_1^{1/2}}(\tau (1,1)^{\top})=\{{\bf 0}_2\}$ when $\tau<3/2^{4/3}\approx 1.19$; and $\prox_{\|\cdot\|_1^{1/2}}(\tau (1,1)^{\top})=\{{\bf 0}_2,(2^{-1/3},2^{-1/3})^{\top}\}$ when $\tau=3/2^{4/3}$. In addition,
 when $\by=(\frac32,y_2)^{\top}$ with $0\le y_2\le \frac13$, then
\(
\prox_{\nu\|\cdot\|_1^{1/2}}(\by)=\{{\bf 0}_2, (1,0)^{\top}\}.
\)
In fact, take $s=1\in[2]$, then $\|\by\|_{(s)}=|y_1|=\frac32\ge\frac{3}{2^{\frac43}}:=\widetilde{t}_{1,1/2}(1)$.  By \eqref{S12ycs}, we arrive at
$$
c_1=\frac{\sqrt{3}}{4\sqrt{\frac32} \cos\big(\frac{1}{3}\arccos(-\frac{3\sqrt{3}}{4}(\frac32)^{-\frac32}) \big)}=\frac{1}{2\sqrt{2}\cos\big(\frac{1}{3}\arccos(-\frac{1}{\sqrt{2}}))}=\frac12.
$$
Clearly, $y_1>c_1$ and $y_2\le c_1$. Hence,
$S_{1/2}(\by,1)=\{\tilde{\by}^{(1)}\}$ with
$\tilde{\by}^{(1)}:=(\frac32-c_1,0)^{\top}=(1,0)^{\top}$  and $J_{\by}(\tilde{\by}^{(1)})=(1+0)^{\frac12}+\frac12((3/2-1)^2+y_2^2)=\frac12y_2^2+\frac98$.
Take $s=2\in[2]$,  $\|\by\|_{(s)}=\|\by\|_1=\frac32+y_2\le \frac32+\frac13=\frac{11}{6}< \frac{3}{2^{\frac23}}=\widetilde{t}_{1,1/2}(2)$, and then  $S_{1/2}(\by,2)=\emptyset$. From \eqref{Spys}, it holds that $S_q(\by)=\{\tilde{\by}^{(1)}\}$. Notice that 
$J_{\by}({\bf 0}_2)=\frac12(\frac94+y_2^2)=\frac{9}{8}+\frac12 y_2^2$. The conclusion holds as $J_{\by}({\bf 0}_2)=J_{\by}(\tilde{\by}^{(1)})$. Furthermore, we compute $\prox_{\|\cdot\|_1^{1/2}}(\tau (3/2,1/2)^{\top})=\{{\bf 0}_2,(1,0)^{\top}\}$.
In summary, $\prox_{\nu\|\cdot\|_1^{1/2}}$ with $\nu=1$ can be verified in Figure \ref{ProxL11over2}.

\subsection{ Characterization of $\prox_{\nu\|\cdot\|_1^{2/3}}$}
  
Similarly, by \eqref{Spy} with $q=2/3$,  we obtain
\begin{equation}\label{S23y}
    S_{2/3}(\by)\!=\!\bigcup_{s\in [n]}\Big\{\max\{\by-c_s{\bf 1}_n ,{\bf 0}_n \}:  \|\by\|_{(s)}\ge \Big(\frac43\frac{2^\frac{3}{4}}{3^{\frac12}}\nu^{\frac34}\Big)s^{\frac34}, y_s\!>\!c_s,  y_{i}\!\le\! c_s, \forall i\!\in\! \{s+1,\ldots,n\}\Big\},
\end{equation}
where  $c_s$ if $\|\by\|_{(s)}\ge \Big(\frac43\frac{2^\frac{3}{4}}{3^{\frac12}}\nu^{\frac34}\Big)s^{\frac34}$ is  by
\begin{equation}\label{S23ycs}
    c_s=\frac{2\nu }{3} \Big(\prox_{\nu s|\cdot|^{2/3}}(\|\by\|_{(s)})\Big)^{-1/3}=\frac{4\nu}{3} \frac{1}{\sqrt{2\alpha}+\sqrt{\frac{2\|\by\|_{(s)}}{\sqrt{2\alpha}}-2\alpha}}
\end{equation}
with $\alpha:=\Big(\frac{\|\by\|_{(s)}^2}{16}+\sqrt{\frac{\|\by\|_{(s)}^4}{256}-\frac{8\nu^3s^3}{729}}\Big)^{1/3}+\Big(\frac{\|\by\|_{(s)}^2}{16}
-\sqrt{\frac{\|\by\|_{(s)}^4}{256}-\frac{8\nu^3s^3}{729}}\Big)^{1/3}$.
  
\begin{corollary}\label{Corollary2over3}
Given $\by\in \mathbb{R}^n_{\downarrow}$. Then
\begin{align}\label{q=2/3Resut}
    \prox_{\nu\|\cdot\|_1^{2/3}}(\by)=
    \left\{\begin{array}{ll}
    \{{\bf 0}_n\}, & \mbox{ if } \|\by\|_1\le 4\big(\frac{2}{9}\big)^{3/4}\nu^{3/4}, \\
    \{\by-c_n{\bf 1}_n\}, & \mbox{ if } \by\succ 2(\frac{2}{3}\nu)^{3/4}{\bf 1}_n,\\
    \arg\min\limits_{\bu\in S_{2/3}(\by)} J_{\by}(\bu), &\mbox{ if }\|\by\|_{\infty}> 2(\frac{2}{3}\nu)^{3/4}\mbox{ and }\by\nsucc 2(\frac{2}{3}\nu)^{3/4}{\bf 1}_n,\\
    \arg\min\limits_{\bu\in\{{\bf 0}_n\}\cup S_{2/3}(\by)} J_{\by}(\bu),&\mbox{ otherwise},
    \end{array}
    \right.
\end{align}
where $S_{2/3}(\by)$ is given by \eqref{S23y} and $c_n$ is defined as in \eqref{S23ycs}.
\end{corollary}
  
\begin{remark}
Hu el. at \cite[Proposition 18 (vi)]{Hu2017} or \cite[Proposition 4.1 (vi)]{Hu2024} proposed the following expression of $\prox_{\nu\|\cdot\|_1^{2/3}}(\by)$ for any $\by\in \mathbb{R}^n$ without proof: 
\begin{equation}\label{wrongL123}
\prox_{\nu\|\cdot\|_1^{2/3}}(\by)=
\left\{\begin{array}{ll}
\{{\bf 0}_n\}, & \mbox{ if }J_{\by}(\tilde{\by}) > J_{\by}({\bf 0}_n),\\
\{{\bf 0}_n,\tilde{\by}\}, & \mbox{ if }J_{\by}(\tilde{\by})=J_{\by}({\bf 0}_n),\\
\{\tilde{\by}\}, & \mbox{ if } J_{\by}(\tilde{\by}) < J_{\by}({\bf 0}_n),
\end{array}
  \right.
\end{equation}
where  $\tilde{\by}:=\by-c_n {\rm sign} (\by)$ with $c_n$ is by \eqref{S23ycs}. Notice that 
$\{\tilde{\by}\}=S_{{\color{black} 2/3 }}(\by,n)$, where $S_{2/3 }(\by,n)$ is by \eqref{Spys}. 
Hence, 
when $S_{2/3}(\by,n)= \arg\min\limits_{\bu\in S_{2/3}(\by)} J_{\by}(\bu)$, 
\eqref{q=2/3Resut} is equivalent to \eqref{wrongL123}. However, the former is more specific; otherwise, the expression \eqref{wrongL123} is incorrect, see the following example.
 \end{remark}

\begin{example}\label{Exmp2} (Counterexample) Fix $\nu=1$ and $n=2$. Let $\by=(1.5,0.7)^{\top}$. Take $s=1\in[2]$, then $\|\by\|_{(s)}=|y_1|=1.5\ge \frac43\frac{2^\frac{3}{4}}{3^{\frac12}}=\widetilde{t}_{1,2/3}(1)\approx 1.2946$. From \eqref{S23ycs}, it follows that $c_1=0.726$. Clearly, $y_1>c_1$ and $y_2\le c_1$. Thus, $S_{2/3}(\by,1)=\{\bu^1\}$ with $\bu^1:=(1.5-c_1,0)^{\top}=(0.774,0)^{\top}$ and $J_{\by}(\bu^1)=(0.774+0)^{2/3}+0.5((1.5-0.774)^2+(0.7-0)^2)=1.3515$.  Take $s=2\in[2]$, then $\|\by\|_{(s)}=\|\by\|_1=2.2>\frac43\frac{2^\frac{3}{2}}{3^{\frac12}}=\widetilde{t}_{1,2/3}(2)\approx 2.1773$. By \eqref{S23ycs}, we obtain $c_2=0.753$. Since $y_2=0.7<c_2=0.753$, $S_{2/3}(\by,2)=\emptyset$. Hence, $S_{2/3}(\by)=\{\bu^1\}$ by \eqref{Spy}. Notice that $J_{\by}({\bf 0}_2)=0.5(1.5^2+0.7^2)=1.37$. It follows that $\prox_{\nu\|\cdot\|_1^{{\color{black} 2/3 }}}(\by)=\{\bu^1\}$ as $J_{\by}(\bu^1)<J_{\by}({\bf 0}_2)$.

However, write $\bu^2:=\by-c_2(1,1)^{\top}=(0.747,-0.053)^{\top}$ and then $J_{\by}(\bu^2)=(0.747+0.053)^{2/3}+0.5((1.5-0.747)^2+(0.7+0.053)^2)\approx 1.4288>J_{\by}({\bf 0}_2)$.  Hence, from  \eqref{wrongL123}, it follows that $\prox_{\nu\|\cdot\|_1^{2/3}}(\by)=\{{\bf 0}_2\}$ , which is obviously incorrect as $J_{\by}(\bu^1)<J_{\by}({\bf 0}_2)$.
\end{example}

\begin{figure}[!htpb]
  \centering
  \begin{subfigure}[b]{0.31\textwidth}
    \centering
    \includegraphics[width=\linewidth]{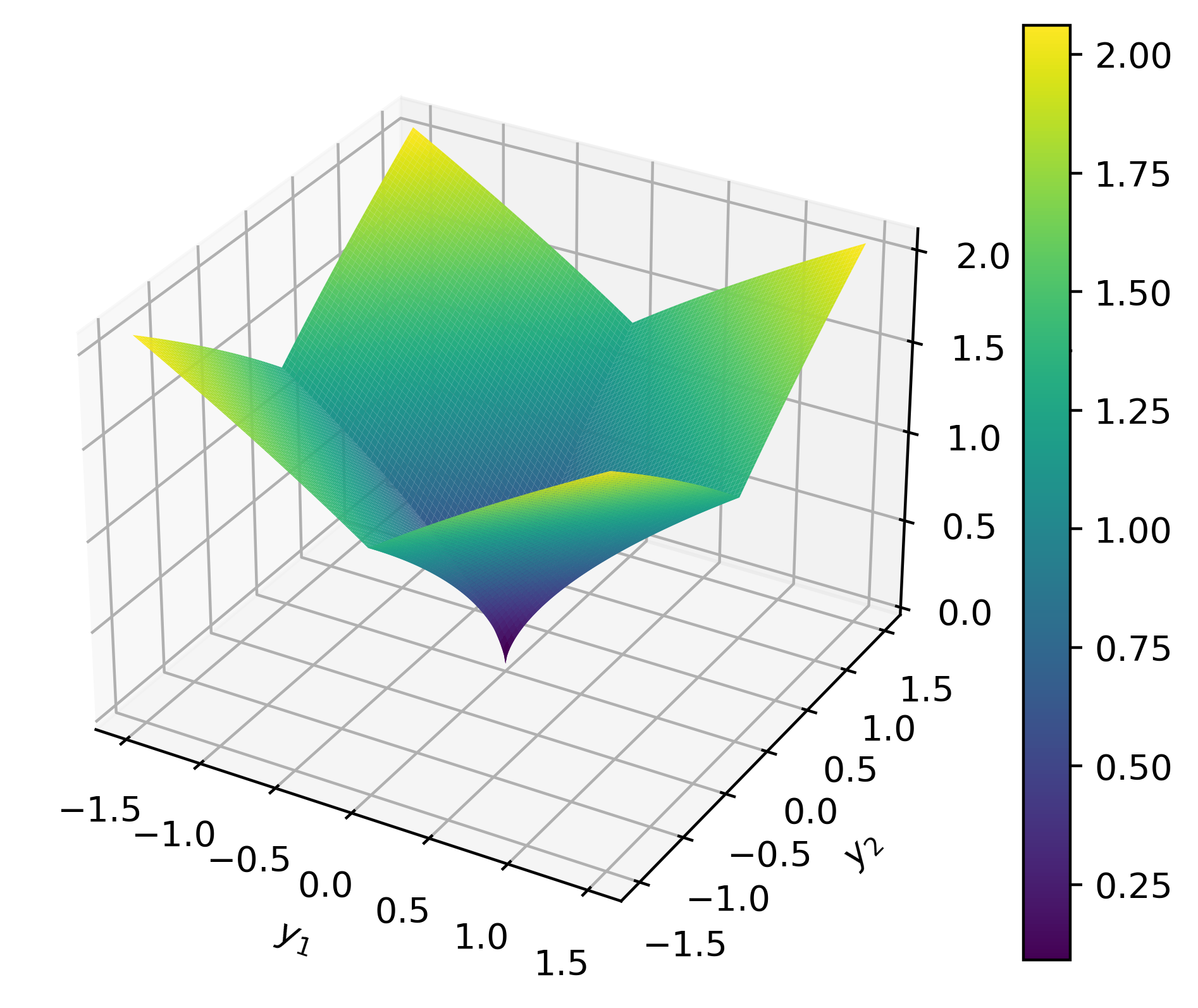}
    \subcaption{Image of $\|\cdot\|_1^{2/3}$ on $\by\in[-2,2]^2$}
    \label{ProxL12over3:sub1}
  \end{subfigure}%
  \hfill
  \begin{subfigure}[b]{0.38\textwidth}
    \centering
    \includegraphics[width=\linewidth]{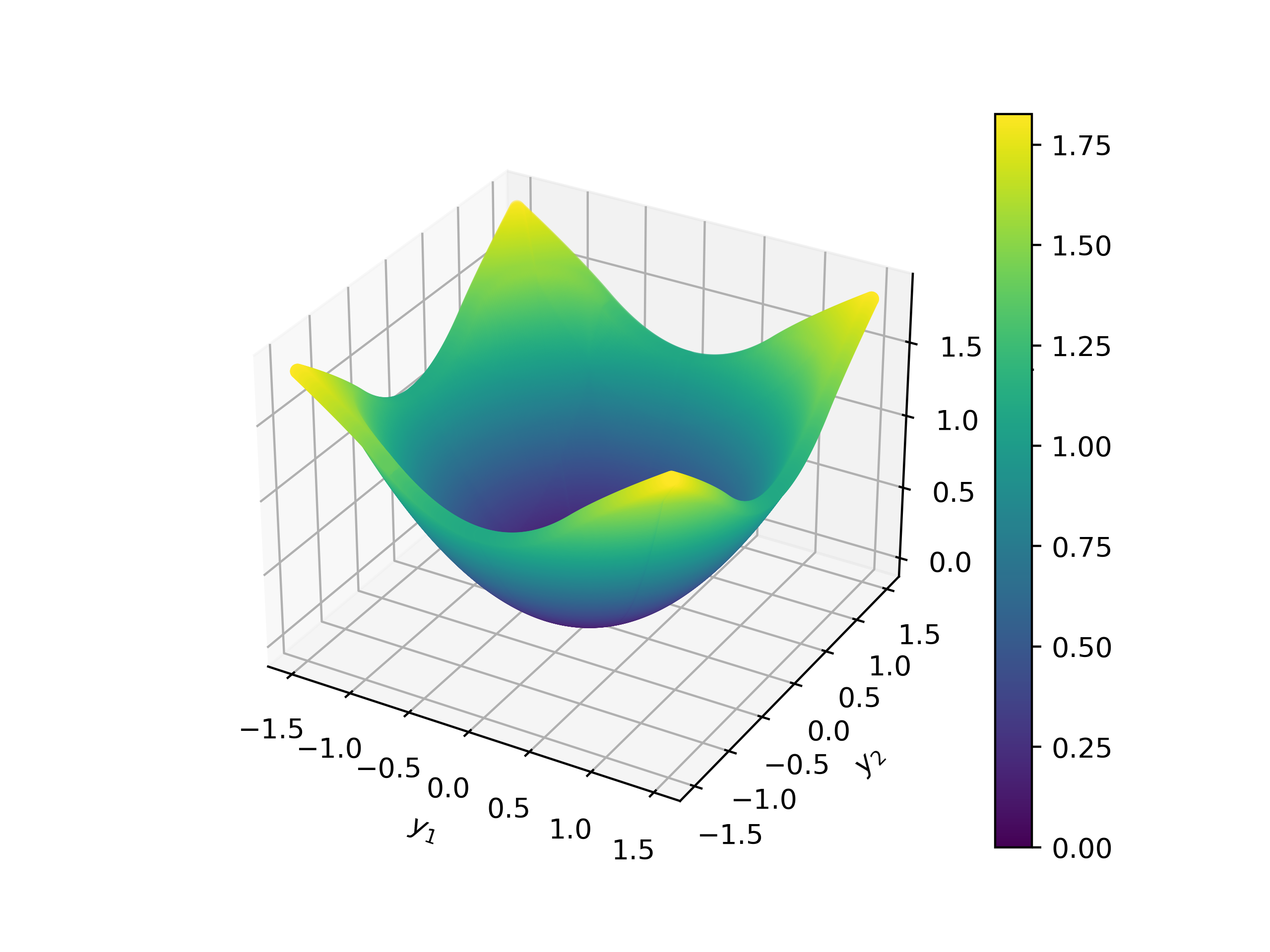}
    \subcaption{Its Moreau envelope}
    \label{ProxL12over3:sub2}
  \end{subfigure}%
  \hfill
  \begin{subfigure}[b]{0.3\textwidth}
    \centering
    \includegraphics[width=\linewidth]{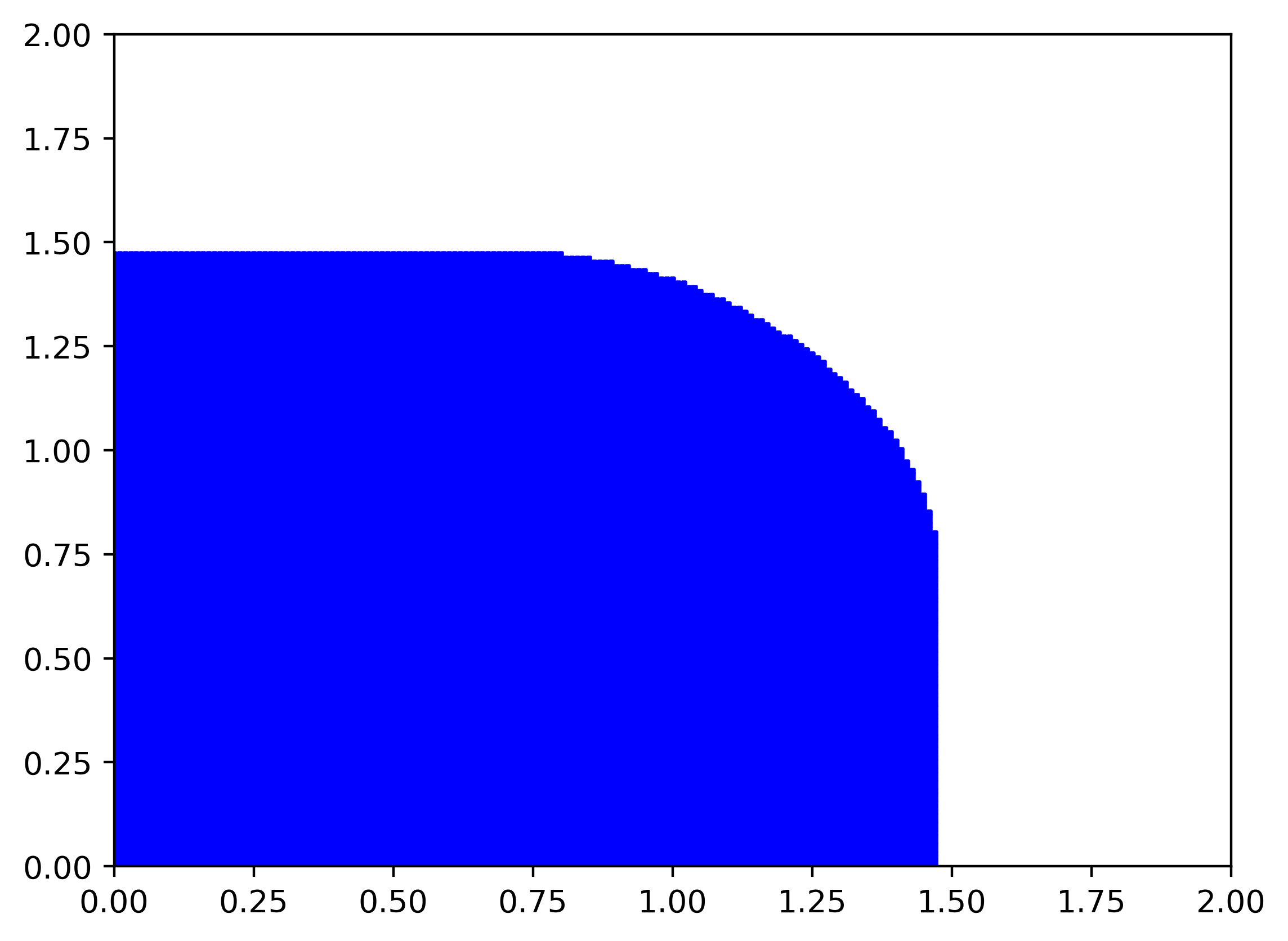}
    \subcaption{$\prox_{\|\cdot\|_1^{2/3}}(\by)\!=\!\{{\bf 0}_2\}$}
    \label{ProxL12over3:sub3}
  \end{subfigure}
  
  \caption{Illustration of $\|\cdot\|_1^{2/3}$, its Moreau envelope and proximal operator.}
  \label{ProxL12over3}
\end{figure}

To understand the behavior of $\prox_{\nu\|\cdot\|_1^{2/3}}$, we make some comments here.  Fix $\nu\!=\!1$.
By Lemma \ref{LpProx} and Lemma \ref{Lemmabasis} , we have  $$\prox_{\|\cdot\|_1^{2/3}}(\tau (1,0)^{\top})=\{{\bf 0}_2\}$$ for any $\tau<2(2/3)^{3/4}\approx 1.4756$; $\prox_{\|\cdot\|_1^{2/3}}(\tau (1,0)^{\top})=\{{\bf 0}_2,((2/3)^{3/4},0)^{\top}\}$ when $\tau=2(2/3)^{3/4}$,  $\prox_{\|\cdot\|_1^{2/3}}(\tau (1,1)^{\top})=\{{\bf 0}_2\}$ for any $\tau<2^{3/2}/3^{3/4}\approx 1.24$,
and $\prox_{\|\cdot\|_1^{2/3}}(\tau (1,1)^{\top})=\{{\bf 0}_2,(2^{1/2}/3^{3/4},2^{1/2}/3^{3/4})^{\top}\}$ when $\tau=2^{3/2}/3^{3/4}$. In addition, 
if $\by\!=\!(2(\frac23)^{3/4},y_2)^{\top}$ with $0\le y_2< \frac43\frac{2^\frac{3}{2}}{3^{\frac12}}\!-\! 2(\frac{2}{3})^{3/4}\approx0.7017$, then 
\(
\prox_{\nu\|\cdot\|_1^{2/3}}(\by)\!=\!\{{\bf 0}_2, ((2/3)^{3/4},0)^{\top}\}
\)
can be obtained via a computation similar to that in Example \ref{Exmp2}. Furthermore, we compute $\prox_{\nu\|\cdot\|_1^{2/3}}((2(2/3)^{3/4},(2/3)^{3/4})^{\top})=\{{\bf 0}_2, ((2/3)^{3/4},0)^{\top}\}$ with $(2/3)^{3/4}\approx 0.7378$.
Hence, the characterization of $\prox_{\nu\|\cdot\|_1^{2/3}}$ with $\nu=1$ can be verified in Figure \ref{ProxL12over3}.

\section{Numerical Experiments}\label{Section4}

The $\ell^q_{1,q}$-based regularization method has been performed for simulated data and real data including gene transcriptional regulation \cite{Hu2017}, group logistic regression \cite{ZhangWei2022}, low-rank matrix completion \cite{Liu2023}, and cell fate conversion \cite{Hu2024}.
The purpose of this section is to implement the numerical experiments on the group sparse vector recovery via the $\ell_{1,q}^{q}$ regularization by using our correct characterization of the proximal operator of the $\ell_{1,q}^{q}$-norm.  
All numerical experiments are implemented in Python 3.11 on a ThinkStation (AMD Ryzen Threadripper Pro 5975WX 32-Cores, 3.6 GHz, 256 GB RAM). 

For comparisons, we also consider group sparse vector recovery via the $\ell_{2,q}^{q}$ regularization. We apply the accelerated proximal gradient algorithm \cite{Beck2009} to solve the $\ell_{p,q}^{q}$ regularization problem \eqref{ProbReg}, which generates a sequence $\{\bx^{(k)}\}\subseteq\bR^l$ via the following iteration:
\begin{equation}\label{PGA}
\left\{\begin{array}{ll}
& \bx^{k}:=\prox_{\lambda\nu\|\cdot\|_{p,q}^q}\Big( \big(\bI_l-\nu \bA^{\top}\bA\big)\by^{k}-\nu \bA^{\top}\bb\Big),\\
& t_{k+1}:=\frac{1+\sqrt{1+4 t_k^2}}{2},\\
& \by^{k+1}:=\bx^{k}+\frac{t_k-1}{t_{k+1}}(\bx^{k}-\bx^{k-1}),
\end{array}\right.
\end{equation}
where $\nu,\lambda>0$, $t_1=1$, the initial point $\by^{1}=\bx^{0}\in\bR^l$ and $\bI_l$ is an $l\times l$ identity matrix. Remember that $\|\bx\|^q_{p,q}\!=\!\sum_{i=1}^r \|\bx_{\mathcal{G}_i}\|_p^q$ with respect to the prescribed partition $\mathcal{G}$. Then, 
$\prox_{\lambda\nu\|\cdot\|_{1,q}^q}$ can be computed by \eqref{MainProx}. It has been proved in \cite[Remark 3.7]{Lin2025} that the proximal operator of the $\ell_{2,q}^{ q }$-norm takes the form  $
 \prox_{\lambda\nu\|\cdot\|^q_{2,q}}(\by)= \prox_{\lambda\nu\|\cdot\|_2^q}(\by_{\mathcal{G}_1})\times \prox_{\lambda\nu\|\cdot\|_2^q}(\by_{\mathcal{G}_2})\times\cdots\times \prox_{\lambda\nu\|\cdot\|_2^q}(\by_{\mathcal{G}_r})
 $ with
$$
\prox_{\lambda\nu\|\cdot\|_2^q}(\by_{\mathcal{G}_i}):=\left\{\begin{array}{ll}
\{{\bf 0}_{|\mathcal{G}_i|}\}, &\mbox{ if } \|\by_{\mathcal{G}_i}\|_2< c_{\lambda\nu,q},\\
\{{\bf 0}_{|\mathcal{G}_i|},\frac{2(1-q)}{2-q}\by_{\mathcal{G}_i}\}, &\mbox{ if } \|\by_{\mathcal{G}_i}\|_2= c_{\lambda\nu,q},\\
\prox_{\lambda\nu\|\by_{\mathcal{G}_i}\|_2^{q-2}|\cdot|^q}(1)\by_{\mathcal{G}_i},& \mbox{ otherwise},
\end{array}i\in[r],\right.
$$
where $c_{\lambda\nu,q}$ is defined in Lemma \ref{LpProx} and $|\mathcal{G}_i|$ denotes the cardinality of $\mathcal{G}_i$. In particular, the proximal operator of the $\ell_{2,q}^{ q }$-norm  has a closed-form expression for $q\in\{1/2,2/3\}$ by Lemma \ref{LpProx}.

Next, we illustrate the performance of the accelerated proximal gradient algorithm \eqref{PGA} among five interesting types of $\ell_{p,q}^{q}$ regularization problem \eqref{ProbReg}, in particular, when $(p,q)=(2,1),(2,2/3),(2,1/2),(1,2/3),(1,1/2)$.

We generate an independent and identically distributed Gaussian matrix $\bA\in\bR^{m\times l}$ with $m=256$ and $l=1024$. All columns of $\bA$ are normalized to unit length with respect to the $\ell_2$-norm. Then, by randomly dividing its components into $r=128$ groups of length $8$ and selecting $k\in[r]$ of them as active groups-whose entries are likewise randomly generated as independent and identically distributed Gaussian—we create a group sparse solution $\bar{\bx}$. The remaining groups are all set to zeros.  The noisy vector $\bb=\bA\bar{\bx}+\sigma*randn(m,1)$, where the parameter $\sigma=0.001$ is the standard deviation of additive Gaussian noise. When the relative error $\|\bx-\bar{\bx}\|_2/\|\bar{\bx}\|_2$ between the true and recovered data is less than 0.5\%, the recovery is considered successful; if not, it is considered failed. For each given inter-group sparsity level $k/r$, we randomly generate $\bA,\bar{\bx},\bb$ as above $100$ times, run the algorithm, and average the $100$ numerical results to illustrate the success rate of different types of $\ell_{p,q}^{ q }$ regularization. For simplicity, we fix $\lambda=1$ and $\nu=0.05$ in the accelerated proximal gradient algorithm \eqref{PGA}.

We first introduce the experiment on signal recovery with only inter-group sparsity. 
Figure \ref{convergence} illustrates the iterative convergence rates of $\ell_{p,q}^{q}$ regularization {\color{black} under four inter-group sparsity levels $k/r$ with $k=1 $, $13$, $19$, and $26$. They correspond to inter-group sparsity levels of approximately 1\%, 10\%, 15\%, and 20\%, respectively. } 
We can see that the proposed algorithm via the $\ell_{1,q}^{q}$ regularization converges faster than the one via $\ell_{2,q}^{ q }$ regularization when the iteration number is relatively small.
It is indicated in Figure \ref{FISTA_success} that the $\ell_{1,1/2}^{1/2 }$- and $\ell_{2,1/2}^{1/2}$-norms can achieve the higher successful recovery rate than other algorithms {\color{black} whenever the inter-group sparsity level greater than $21/128\approx 16.41\%$.} 

{\color{black} Next, we will discuss the experiment on signal recovery with both inter-group and intra-group sparsity.
In this setting, we always fix the intra-group sparsity level to 50\%, such that the proportion of non-zero elements within active groups is 50\%.
Figures \ref{intra_convergence_0.001} and \ref{intra_convergence_0.01} depict the iterative convergence rates for four inter-group sparsity levels with noise levels $\sigma = 0.001$ and $\sigma = 0.01$, respectively. 
Figures \ref{FISTA_success_0.001} and \ref{FISTA_success_0.01} show the successful recovery rate for different inter-group sparsity levels with noise levels $\sigma = 0.001$ and $\sigma = 0.01$, respectively. As the number $k$ of activated groups increases, the $\ell_{1,q}^q$ regularization method demonstrates a smaller recovery error and attains a significantly higher recovery success rate in comparison to the $\ell_{2,q}^q$ regularization method for inter-group and intra-group sparse vector recovery. In terms of overall performance, the $\ell_{1,2/3}^{2/3}$ regularization method attains the highest recovery success rate.
Figures \ref{FISTA_success} and \ref{FISTA_success_0.001}, both with the same noise level of $\sigma=0.001$, illustrate that the success rate of the $\ell_{1,q}^q$ regularization method is greatly enhanced when intra-group sparsity is incorporated.

\begin{figure}[!htpb]
\vspace{-2.3cm}
  \centering
  \begin{subfigure}[b]{0.49\linewidth}
    \centering
    \includegraphics[width=\linewidth]{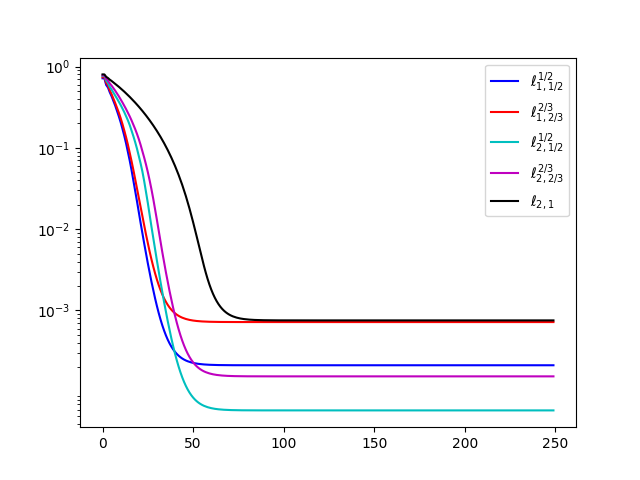}
    \subcaption{$k = 1$}
    \label{convergence:spa_1}
  \end{subfigure}%
  \begin{subfigure}[b]{0.49\linewidth}
    \centering
    \includegraphics[width=\linewidth]{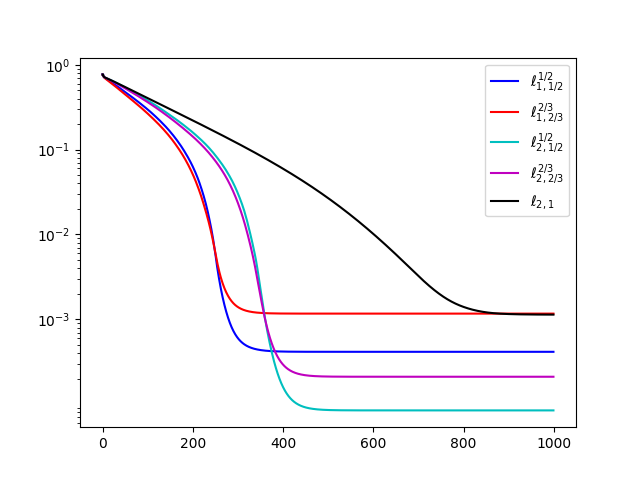}
    \subcaption{$k = 13$}
    \label{convergence:spa_13}
  \end{subfigure}
  \begin{subfigure}[b]{0.49\linewidth}
    \centering
    \includegraphics[width=\linewidth]{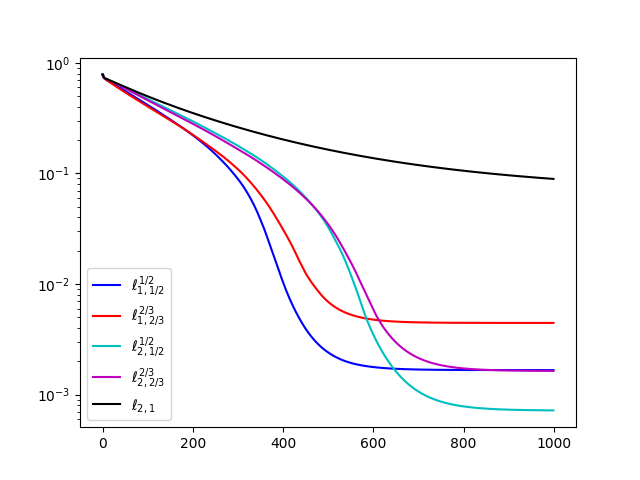}
    \subcaption{$k = 19$}
    \label{convergence:spa_19}
  \end{subfigure}
  \begin{subfigure}[b]{0.49\linewidth}
    \centering
    \includegraphics[width=\linewidth]{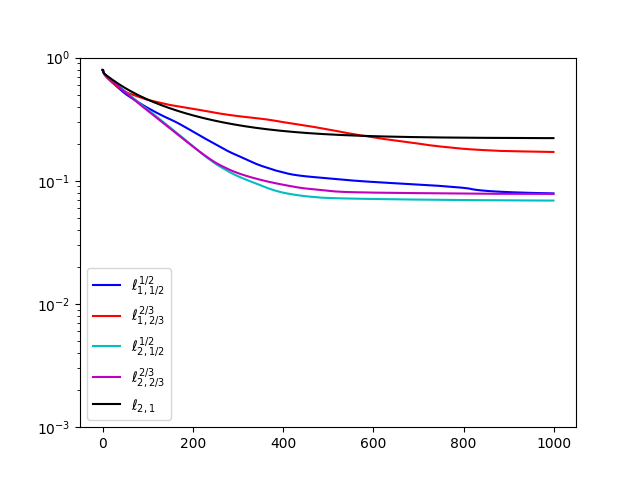}
    \subcaption{$k = 26$}
    \label{convergence:spa_26}
  \end{subfigure}
  
  \caption{The relative error of $\ell_{p,q}^{{\color{black} q }}$ regularization methods versus the number of iterations for $\sigma=0.001$ and four inter-group sparsity levels $k/r$ with $k\in\{1,13,19,26\}$.}
  \label{convergence}
\end{figure}

\begin{figure}[!htpb]
  \centering
   \includegraphics[width=3in]{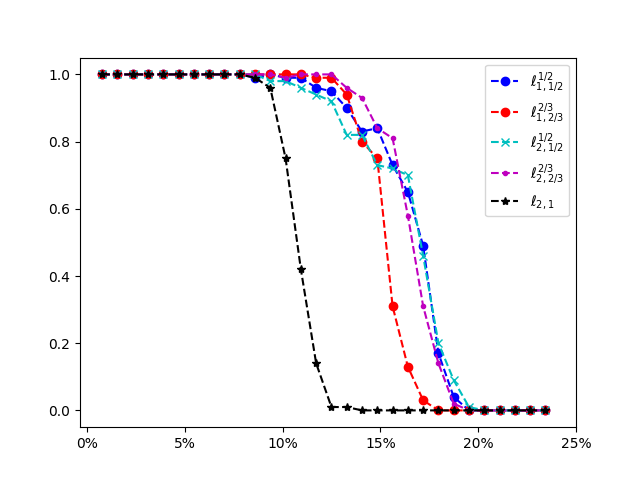}
  \caption{Success rate of $\ell_{p,q}$ regularization methods for $\sigma=0.001$ and different inter-group sparsity levels. }\label{FISTA_success}
\end{figure}

\begin{figure}[!ht]
    \vspace{-2.3cm}
  \centering
  \begin{subfigure}[b]{0.49\linewidth}
    \centering
    \includegraphics[width=\linewidth]{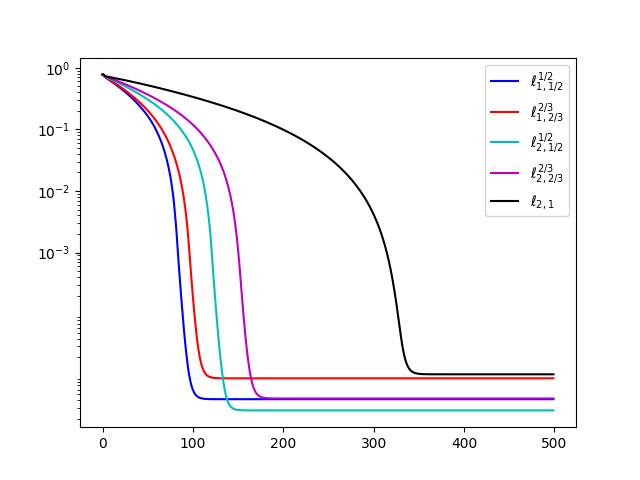}
    \subcaption{$k = 1$}
    \label{intra_convergence_0.001:intra_a}
  \end{subfigure}
  \begin{subfigure}[b]{0.49\linewidth}
    \centering
    \includegraphics[width=\linewidth]{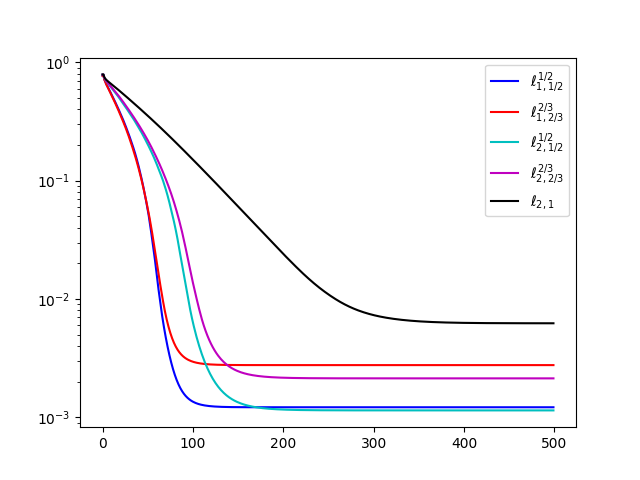}
    \subcaption{$k = 13$}
    \label{intra_convergence_0.001:intra_b}
  \end{subfigure}
  \begin{subfigure}[b]{0.49\linewidth}
    \centering
    \includegraphics[width=\linewidth]{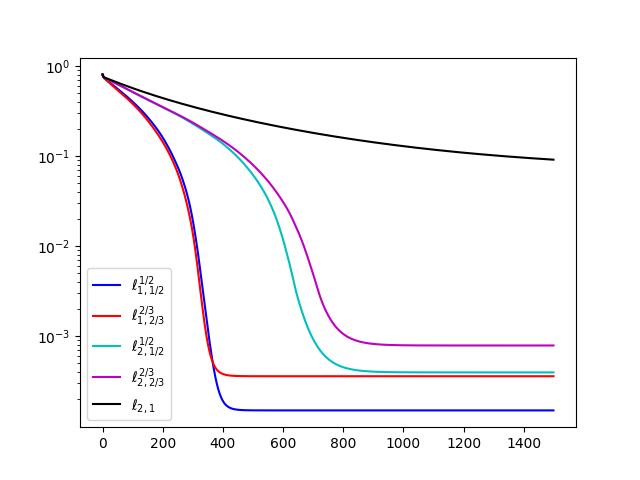}
    \subcaption{$k = 19$}
    \label{intra_convergence_0.001:intra_c}
  \end{subfigure}
  \begin{subfigure}[b]{0.49\linewidth}
    \centering
    \includegraphics[width=\linewidth]{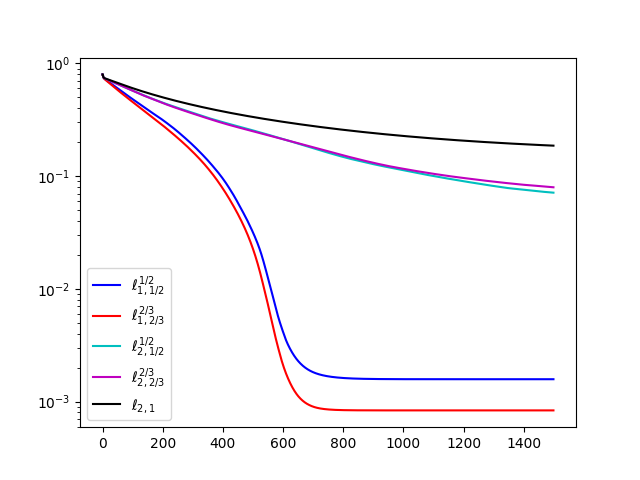}
    \subcaption{$k = 26$}
    \label{intra_convergence_0.001:intra_d}
  \end{subfigure}
  
  \caption{{\color{black}The relative error of of $\ell_{p,q}^{{\color{black} q }}$ regularization methods versus the number of iterations for $\sigma = 0.001$, 50\% intra-group sparsity, and four inter-group sparsity levels $k/r$ with $k\in\{1,13,19,26\}$.}}
  \label{intra_convergence_0.001}
\end{figure}

\begin{figure}[!ht]
  \centering
   \includegraphics[width=3in]{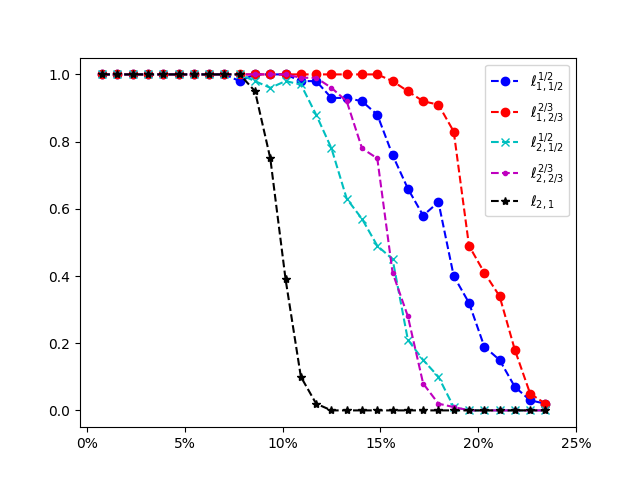}
  \caption{{\color{black}Success rate of $\ell_{p,q}^{{\color{black} q }}$ regularization methods for $\sigma=0.001$ and different sparsity levels. }}\label{FISTA_success_0.001}
\end{figure}

\begin{figure}[!ht]
    \vspace{-2.3cm}
  \centering
  \begin{subfigure}[b]{0.49\linewidth}
    \centering
    \includegraphics[width=\linewidth]{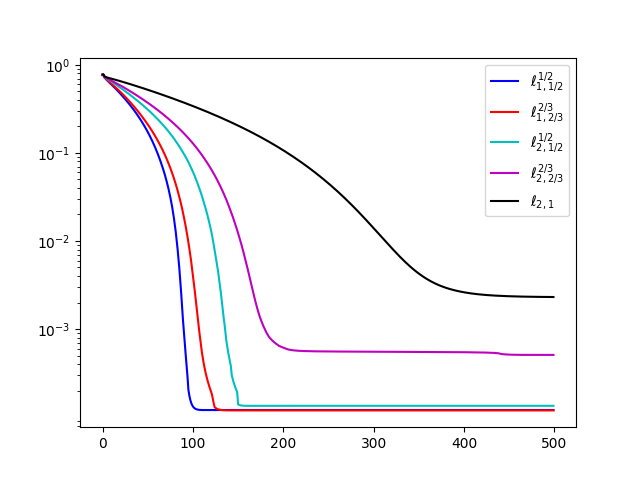}
    \subcaption{$k = 1$}
    \label{intra_convergence_0.01:intra_a}
  \end{subfigure}
  \begin{subfigure}[b]{0.49\linewidth}
    \centering
    \includegraphics[width=\linewidth]{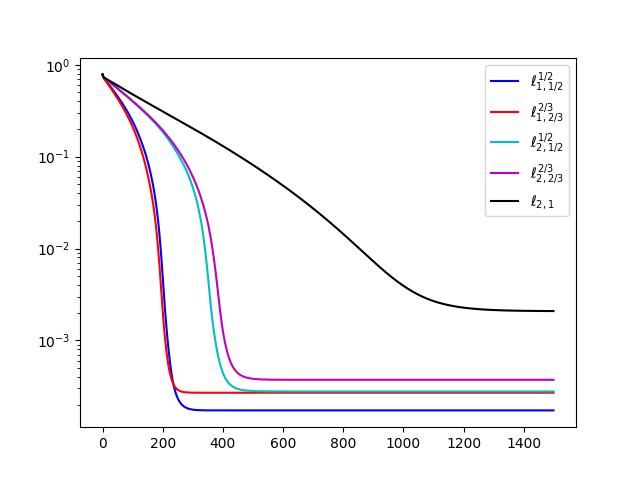}
    \subcaption{$k = 13$}
    \label{intra_convergence_0.01:intra_b}
  \end{subfigure}
  \begin{subfigure}[b]{0.49\linewidth}
    \centering
    \includegraphics[width=\linewidth]{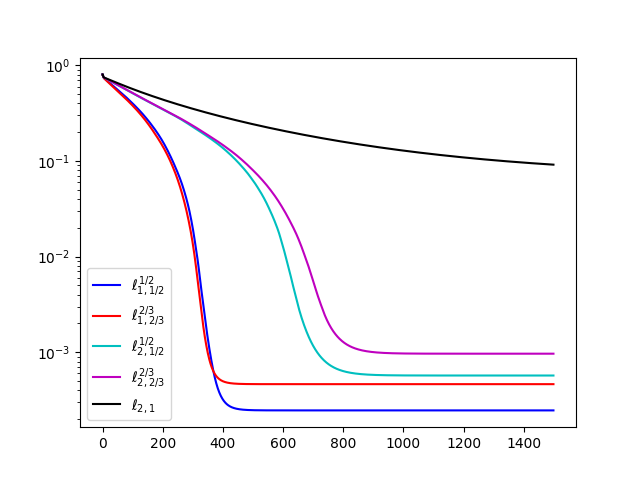}
    \subcaption{$k = 19$}
    \label{intra_convergence_0.01:intra_c}
  \end{subfigure}
  \begin{subfigure}[b]{0.49\linewidth}
    \centering
    \includegraphics[width=\linewidth]{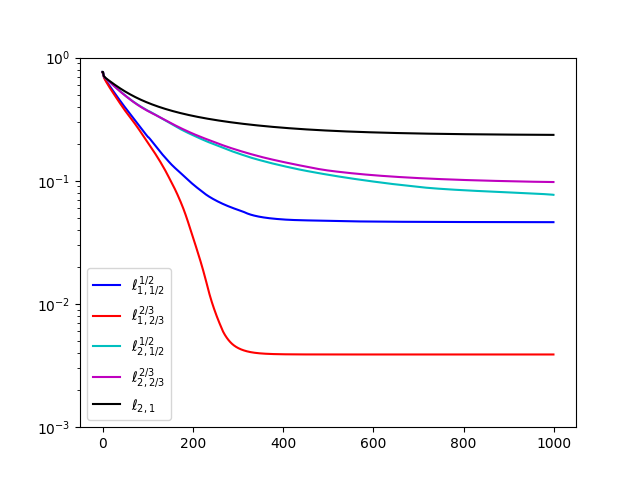}
    \subcaption{$k = 26$}
    \label{intra_convergence_0.01:intra_d}
  \end{subfigure}
  
  \caption{{\color{black}The relative error of $\ell_{p,q}^q$ regularization methods versus the number of iterations for $\sigma = 0.01$, 50\% intra-group sparsity, and four inter-group sparsity levels $k/r$ with $k\in\{1,13,19,26\}$.}}
  \label{intra_convergence_0.01}
\end{figure}

\begin{figure}[!ht]
  \centering
   \includegraphics[width=3in]{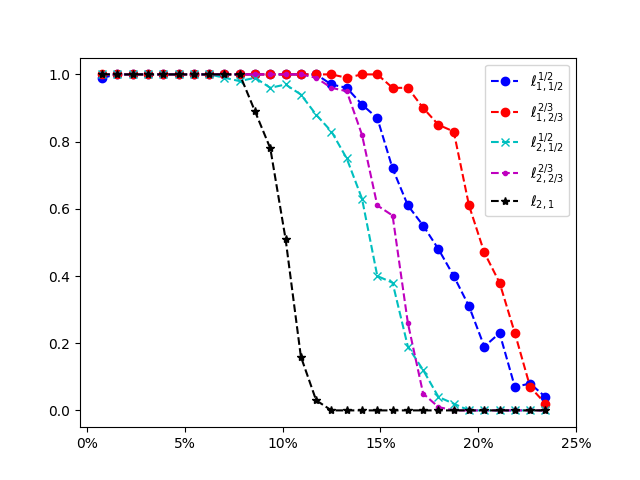}
  \caption{{\color{black}Success rate of $\ell_{p,q}^q$ regularization methods for $\sigma=0.01$ and different sparsity levels. }}\label{FISTA_success_0.01}
\end{figure}
}

\section{Conclusions and remarks}

In this note, we derive a characterization in Theorem \ref{MainTh} for the proximal operator of the $\ell^{ q}_{1,q}$-norm with $0\!<\!q\!<\!1$ by means of the proximal operator of  the $\ell^q_q$-norm. In particular, the proximal operators of the $\ell_{1,1/2}^{ 1/2 }$-norm and the $\ell_{1,2/3}^{ 2/3 }$-norm are presented in Corollary \ref{Corollary1over2} and Corollary \ref{Corollary2over3}, respectively. {\color{black}Numerical experiments in Figures \ref{FISTA_success_0.001} and \ref{FISTA_success_0.01} have shown that the $\ell_{1,q}^q$ regularization method is highly effective in addressing both inter-group and intra-group sparse issues.}
This work will help develop efficient proximal algorithms for group sparse optimization with the $\ell_{1,q}^{ q }$ regularization in various areas including gene expression analysis, signal and image processing, low-rank matrix completion, multiple kernel learning, neural networks, etc.



\bibliographystyle{elsarticle-num}
\bibliography{references}

\begin{thebibliography}{10}
\expandafter\ifx\csname url\endcsname\relax
  \def\url#1{\texttt{#1}}\fi
\expandafter\ifx\csname urlprefix\endcsname\relax\def\urlprefix{URL }\fi
\expandafter\ifx\csname href\endcsname\relax
  \def\href#1#2{#2} \def\path#1{#1}\fi

\bibitem{Meier2008}
L.~Meier, S.~Van De~Geer, P.~B{\"u}hlmann, The group lasso for logistic regression, J. R. Stat. Soc. Ser. B Stat. Methodol. 70~(1) (2008) 53--71.

\bibitem{Hu2024}
Y.~Hu, X.~Hu, C.~K.~W. Yu, J.~Qin, Joint sparse optimization: lower-order regularization method and application in cell fate conversion, Inverse Probl. 40~(9) (2024) 095003.

\bibitem{Sun2014}
Y.~Sun, Q.~Liu, J.~Tang, D.~Tao, Learning discriminative dictionary for group sparse representation, IEEE Trans. Image Process. 23~(9) (2014) 3816--3828.

\bibitem{Wang2017}
Q.~Wang, Q.~Gao, X.~Gao, F.~Nie, $\ell_{2,p}$-norm based {PCA} for image recognition, IEEE Trans. Image Process. 27~(3) (2017) 1336--1346.

\bibitem{Liu2023}
Q.~Liu, Q.~Jiang, J.~Zhang, B.~Jiang, Z.~Liu, An improved low-rank matrix fitting method based on weighted ${L}_{1,p}$ norm minimization for matrix completion, Int. J. Pattern Recognit Artif Intell. 37~(04) (2023) 2350007.

\bibitem{ZhangLi2024}
X.~Zhang, X.~Li, C.~Zhang, Matrix optimization problem involving group sparsity and nonnegativity constraints, J. Optim. Theory App. 201~(1) (2024) 130--176.

\bibitem{LiuJi2009}
J.~Liu, S.~Ji, J.~Ye, Multi-task feature learning via efficient $\ell_{2,1}$-norm minimization, in: Proceedings of the 25th Conference on Uncertainty in Artificial Intelligence, UAI 2009, AUAI Press, 2009, pp. 339--348.

\bibitem{Bach2008}
F.~R. Bach, Consistency of the group lasso and multiple kernel learning, J. Mach. Learn. Res. 9~(6) (2008) 1179--1225.

\bibitem{Wen2016}
W.~Wen, C.~Wu, Y.~Wang, Y.~Chen, H.~Li, Learning structured sparsity in deep neural networks, Advances in Neural Information Processing Systems 29 (2016).

\bibitem{Baldassarre2016}
L.~Baldassarre, N.~Bhan, V.~Cevher, A.~Kyrillidis, S.~Satpathi, Group-sparse model selection: Hardness and relaxations, IEEE Trans. Inf. Theory 62~(11) (2016) 6508--6534.

\bibitem{Eldar2009}
Y.~C. Eldar, M.~Mishali, Robust recovery of signals from a structured union of subspaces, IEEE Trans. Inf. Theory 55~(11) (2009) 5302--5316.

\bibitem{Eldar2010}
Y.~C. Eldar, P.~Kuppinger, H.~Bolcskei, Block-sparse signals: Uncertainty relations and efficient recovery, IEEE Trans. Signal Process. 58~(6) (2010) 3042--3054.

\bibitem{Beck2019}
A.~Beck, N.~Hallak, Optimization problems involving group sparsity terms, Math. Program. 178 (2019) 39--67.

\bibitem{Duarte2009}
M.~F. Duarte, V.~Cevher, R.~G. Baraniuk, Model-based compressive sensing for signal ensembles, in: 2009 47th Annual Allerton Conference on Communication, Control, and Computing, 2009, pp. 244--250.

\bibitem{Hu2017}
Y.~Hu, C.~Li, K.~Meng, J.~Qin, X.~Yang, Group sparse optimization via $\ell_{p,q}$ regularization, J. Mach. Learn. Res. 18~(30) (2017) 1--52.

\bibitem{Pan2021}
L.~Pan, X.~Chen, Group sparse optimization for images recovery using capped folded concave functions, SIAM J. Imaging Sci. 14~(1) (2021) 1--25.

\bibitem{Yuan2006}
M.~Yuan, Y.~Lin, Model selection and estimation in regression with grouped variables, J. R. Stat. Soc. Ser. B Stat. Methodol. 68~(1) (2006) 49--67.

\bibitem{Zhang2022}
X.~Zhang, J.~Zheng, D.~Wang, G.~Tang, Z.~Zhou, Z.~Lin, Structured sparsity optimization with non-convex surrogates of $\ell_{2,0}$-norm: A unified algorithmic framework, IEEE Trans. Pattern Anal. Mach. Intell. 45~(5) (2022) 6386--6402.

\bibitem{Feng2020}
X.~Feng, S.~Yan, C.~Wu, The $\ell_{2,q}$ regularized group sparse optimization: lower bound theory, recovery bound and algorithms, Appl. Comput. Harmon. Anal. 49~(2) (2020) 381--414.

\bibitem{Chen2023}
X.~Chen, L.~Pan, N.~Xiu, Solution sets of three sparse optimization problems for multivariate regression, J. Global Optim. (2023) 1--25.

\bibitem{Beck2017}
A.~Beck, First-order methods in optimization, SIAM, 2017.

\bibitem{Parikh2014}
N.~Parikh, S.~Boyd, et~al., Proximal algorithms, Found. Trends Optim. 1~(3) (2014) 127--239.

\bibitem{Blumensath2008}
T.~Blumensath, M.~E. Davies, Iterative thresholding for sparse approximations, J. Fourier Anal. Appl. 14 (2008) 629--654.

\bibitem{Beck2009}
A.~Beck, M.~Teboulle, A fast iterative shrinkage-thresholding algorithm for linear inverse problems, SIAM J. Imaging Sci. 2~(1) (2009) 183--202.

\bibitem{Cao2013}
W.~Cao, J.~Sun, Z.~Xu, Fast image deconvolution using closed-form thresholding formulas of $\ell_q$ $(q= 1/2,2/3)$ regularization, J. Visual Commun. Image Represent. 24~(1) (2013) 31--41.

\bibitem{Chen2016}
F.~Chen, L.~Shen, B.~W. Suter, Computing the proximity operator of the $\ell_p$ norm with $0< p< 1$, IET Signal Process. 10~(5) (2016) 557--565.

\bibitem{ZhangTong2010}
T.~Zhang, Analysis of multi-stage convex relaxation for sparse regularization, J. Mach. Learn. Res. 11~(3) (2010).

\bibitem{Zhang2017}
S.~Zhang, J.~Xin, Minimization of transformed $\ell_1$ penalty: Closed form representation and iterative thresholding algorithms, Commun. Math. Sci. 15~(2) (2017) 511--537.

\bibitem{Fan2001}
J.~Fan, R.~Li, Variable selection via nonconcave penalized likelihood and its oracle properties, J. Amer. Statist. Assoc. 96~(456) (2001) 1348--1360.

\bibitem{Zhang2010}
C.-H. Zhang, Nearly unbiased variable selection under minimax concave penalty, Ann. Statist. 38~(2) (2010) 894--942.

\bibitem{Liu2024}
Y.~Liu, Y.~Zhou, R.~Lin, The proximal operator of the piece-wise exponential function, IEEE Signal Process. Lett. 31 (2024) 894--898.

\bibitem{Prater2022}
A.~Prater-Bennette, L.~Shen, E.~E. Tripp, The proximity operator of the log-sum penalty, J. Sci. Comput. 93~(67) (2022) 1--34.

\bibitem{Lou2018}
Y.~Lou, M.~Yan, Fast ${L}_1-{L}_2$ minimization via a proximal operator, J. Sci. Comput. 74~(2) (2018) 767--785.

\bibitem{Tao2022}
M.~Tao, Minimization of ${L}_1$ over ${L}_2$ for sparse signal recovery with convergence guarantee, SIAM J. Sci. Comput. 44~(2) (2022) A770--A797.

\bibitem{Jia2024}
J.~Jia, A.~Prater-Bennette, L.~Shen, Computing proximity operators of scale and signed permutation invariant functions, arXiv:2404.00713 (2024).

\bibitem{Prater2023}
A.~Prater-Bennette, L.~Shen, E.~E. Tripp, A constructive approach for computing the proximity operator of the $p$-th power of the $\ell_1$ norm, Appl. Comput. Harmon. Anal. 67 (2023) 101572.

\bibitem{Lin2025}
R.~Lin, S.~Li, Z.~Li, Y.~Liu, Convergence analysis of iteratively reweighted $\ell_1$ algorithms for computing the proximal operator of $\ell_p$-norm, Math. Methods Appl. Sci. 48~(1) (2025) 386--396.

\bibitem{Zwillinger2018}
D.~Zwillinger, CRC Standard Mathematical Tables and Formulas, 33rd Edition, Chapman and Hall/CRC, 2018.

\bibitem{Liu2024b}
Y.~Liu, R.~Lin, A bisection method for computing the proximal operator of the $\ell_p$-norm for any $0<p<1$ with application to {S}chatten $p$-norms, J. Comput. Appl. Math. 447 (2024) 115897.

\bibitem{RW09}
R.~T. Rockafellar, R.~J.-B. Wets, Variational analysis, Vol. 317, Springer Science \& Business Media, 2009.

\bibitem{ZhangWei2022}
Y.~Zhang, C.~Wei, X.~Liu, Group logistic regression models with $\ell_{p,q}$ regularization, Mathematics 10~(13) (2022) 2227.

\end{thebibliography}

\end{document}